\documentclass[reqno]{amsart}

\usepackage{amsfonts,amsmath,amssymb,delarray}
\usepackage{enumerate}
\newcommand{\beq}{\begin{equation}}
\newcommand{\eeq}{\end{equation}}
\newcommand{\be}{\begin{equation*}}
\newcommand{\ee}{\end{equation*}}
\newcommand{\bs}{\begin{split}}
\newcommand{\esplit}{\end{split}}
\newcommand{\begincal}{\begin{eqnarray}}
\newcommand{\fincal}{\end{eqnarray}}

\newtheorem{thm}{Theorem}[section]
\newtheorem{lemma}{Lemma}[section]
\newtheorem{cor}{Corollary}[section]
\newtheorem{prop}{Proposition}[section]
\newtheorem{defi}{Definition}[section]

\newcommand{\ds}{\displaystyle}
\newcommand{\eps}{\varepsilon}

\newcommand{\R}{{\mathbb R}}

\newcommand{\C}{{\mathbb C}}
\newcommand{\N}{{\mathbb N}}

\begin{document}
 
\title[Surfaces with large constant mean curvature]{Concentration of $CMC$ surfaces in a $3$-manifold}
\author{Paul Laurain}
\thanks{Paul Laurain UMPA-ENS Lyon 46 allée d'Italie 69364 Lyon Cedex 07.\\ paul.laurain@umpa.ens-lyon.fr}

\begin{abstract} 
We prove that simply connected $H$-surfaces with small diameter  in a $3$-manifold necessarily concentrate at a critical point of the scalar curvature.
\end{abstract}

\maketitle

\section*{Introduction}
Let $(N,g)$ be a compact  oriented Riemannian manifold. The aim of this article is to understand the behaviour of a sequence of surfaces $\Sigma_H\subset N$ with constant mean curvature $H$, refereed to as $H$-surfaces, when $H\rightarrow +\infty$.\\

These $H$-surfaces naturally appear as boundaries of isoperimetric domains. Their existence  is given by geometric measure theory \cite{Morgan}, but we have no information about their topology or their location in the considered manifold except for some special manifolds like space forms where we have a classification of compact embedded $H$-surfaces (this is an extension of Aleksandrov theorem \cite{Aleksandrov56}, see for instance \cite{MontielRos}).\\

It would be too ambitious for now to hope for a classification of these $H$-surfaces in a general compact manifold. However, in the particular case of minimal surfaces ($H = 0$), a rough classification can be obtained thanks to the works of Colding, Meeks, Minicozzi, Ros and Rosenberg and others. We will find an overview on this subject in the collective book edited by Hoffman \cite{Hoffman} and the papers of Colding and Minicozzi \cite{CM1}, \cite{CM2}, \cite{CM3} and  \cite{CM4}. This area of research is still very active motivated by its close links with the topology of $3$-manifolds.

In order to begin the description of the moduli space of $H$-surfaces, we look to the case of surfaces with small diameter (or large mean curvature). Up to perform dilation of the ambient space, we can normalize the mean curvature of these surfaces to be $1$ and the ambient space becomes quasi-Euclidean. In this setting, an idea to obtain explicit examples of constant mean curvature surfaces is to pertub the constant mean curvature surfaces of the Euclidean space (i.e. round spheres but also connected sums of spheres and Delaunay surfaces)  in order to get surfaces with constant mean curvature in our quasi-Euclidean space. This idea has been very successful and has led to  many examples, see Ye \cite{Ye91}, Butscher \cite{Butscher}, Butscher-Mazzeo \cite{ButscherMazzeo}, Pacard \cite{Pacard} and Pacard-Xu \cite{PacardXu09}. But each of these constructions requires a condition on the geometry of the manifold at the point of concentration. A natural question then is the question of the necessity of this geometric condition. In fact if we were able to show that these conditions are necessary, we would have a clearer picture of the moduli space, at least for surfaces of small diameter. A first answer was given by Druet \cite{Druet02} in the case of isoperimetric domains. By proving an optimal isoperimetric inequality for domains of small volumes, he shows that these domains concentrate necessarily  at a point of maximum scalar curvature. This result, together with the examples mentioned above, leads naturally to the following question, already mentioned in \cite{PacardXu09}~: if, for any $\rho>0$, the ball $B\left(p,\rho\right)$ contains a constant mean curvature surface, is it true that $p$ has to be a critical point of the scalar curvature~? All the examples mentioned above are constructed in a neighbourhood of a critical point of the scalar curvature (with various nondegeneracy assumptions). 

For isoperimetric domains, the topology and geometry of the domains become simple as the volume goes to $0$ (that is as the constant mean curvature goes to $+\infty$). Indeed, they asymptotically become round spheres (see \cite{JohnsonMorgan00, Druet02, Nardulli}). Of course, without this minimizing property of isoperimetric domains, the geometry of constant mean curvature surfaces becomes more intricate. Notably, even in the embedded case, $\Sigma_H$ could be  a connected sum of Delaunay surfaces and an arbitrary number of almost round spheres. Indeed, Pacard and Malchiodi (see \cite{Pacard}) have constructed sequences of $H$-surfaces which are perturbations of two small geodesic spheres connected as a Delaunay surface. Another problem is the topology of the surface which is {\it a priori} unknown, even if the ambient space is Euclidean as shown by Wente's tori and then by Kapouleas's surfaces, see \cite{Wente86} and \cite{Kapouleas90}. In order to generalize the result of Druet \cite{Druet02}, we consider sequences of $H$-surfaces $\Sigma_H$ which are embedded spheres with  bounded area and small diameter. The assumptions are precisely the following :
\beq
\tag{H}
\label{H}
\begin{array}{cc} \begin{cases}
 \delta  (\Sigma_{H})=o(1) \\
 A(\Sigma_{H}) = O\left(\frac{1}{H^2}\right)
\end{cases} & \hbox{ as } H\rightarrow +\infty,  \end{array}
\eeq 
where $\delta  (\Sigma_{H})$  and $A(\Sigma_{H})$ denote respectively the extrinsic diameter and the area of $\Sigma_H$. Here, the area is computed with respect to the induced metric. Under these assumptions, we are able to locate the possible places of concentration of these sequences~:

\begin{thm}
\label{thmain}
Let $(N,g)$ be a smooth compact $3$-Riemannian manifold and $\Sigma_{H}\subset N$ be a sequence of embedded spheres with constant mean curvature $H$ which satisfies assumptions  (\ref{H}). Then, $\Sigma_H$ converges uniformly to a critical point of the scalar curvature.
\end{thm}

We can rephrase this theorem as follows~: choose any function $\delta(H)$ such that $\delta(H)\to 0$ as $H\to +\infty$. Then for any $\rho>0$, there exists some $H_0>0$ such that any embedded topological sphere with constant mean curvature $H\ge H_0$, diameter $\delta\le \delta(H)$ and area $A\le \rho^{-1}H^{-2}$, has to be in a ball $B\left(p,\rho\right)$ where $p$ is a critical point of the scalar curvature of $(N,g)$. 

Conversely, if $p$ is a nondegenerate critical point of the scalar curvature, then there are such embedded spheres in any ball $B\left(p,\rho\right)$ (see Ye \cite{Ye91}). 

Moreover, as will be seen from the proof of the theorem, we get a precise asymptotic description of the surfaces $\Sigma_H$ as $H\to +\infty$~: roughly speaking, they look like a connected sum of spheres. 

This theorem thus provides a beginning of classification of high constant mean curvature surfaces in $3$-dimensional Riemannian manifolds. 

\medskip Note that one could ask the same question for curves in $2$-manifolds. And the answer is simpler than in dimension $3$~: curves with high constant geodesic curvature and small diameter converge to some critical point of the Gauss curvature. This was proved by Sun \cite{Sun08}. The difference between curves and surfaces is that, for curves, one has to analyze solutions of some ODE while, here, we have to deal with solutions of some system of elliptic PDEs. We also note that a similar theorem has been proved by the author concerning small constant mean curvature surfaces with boundary in a euclidean domain, see \cite{Laurain10bis}.

\medskip The rest of the paper is devoted to the proof of theorem \ref{thmain} and is organized as follows. First, in section 1, we compute the equation satisfied by our $H$-surfaces  in a general $3$-manifold and  we recall the classification of solutions of the limit equation (i.e. when the ambient metric becomes flat) obtained by Brezis and Coron \cite{BC2}. In section 2, we set up our proof by reformulating the problem in the framework of some blow-up analysis for a sequence of solutions of perturbed $H$-systems. These systems are systems of elliptic PDEs, critical form the point of view of Sobolev embeddings, but which enjoy some nice compactness by compensation properties (see Rivi\`ere \cite{Riviere} for a nice and clear explanation of these phenomena). However, the perturbation due to the presence of some Riemannian metric, instead of the Euclidean one, breaks most of these properties. In section 3, we start the blow-up analysis by showing that our sequence of solutions decomposes asymptotically into a sum of spheres. This is a generalization of the classical result of Brezis and Coron \cite{BC2} in our setting. Then comes the key point of the proof~: we need to estimate precisely (and in a pointwise way) the error between our solutions and this sum of parametrizations of spheres. Roughly speaking, we have to upgrade the theory of Brezis-Coron which took place in the energy space into a pointwise theory, following the general scheme developed for Yamabe type equations by Hebey et al., see e.g. \cite{HebeyVaugon}, \cite{DruetHebeyRobert}, \cite{Hebey06}, \cite{DruetHebey08}. This is done in two steps. We first use an estimate obtained thanks to the Green formula and a classification of decreasing solutions  of  the linearized equation, see section \ref{lin}. It  remains to control the interaction between the bubbles, which is postponed at the end of proof. Finally, this estimate is in section \ref{s4} to conclude. The proof is rather technical and the reader can start by assuming that there is juste one bubble, like in the construction of Ye. In this case, one can ignore section \ref{s7}. In the general case, when there are several bubbles, we must also get a good control on the interaction between bubbles.
 
 \medskip {\bf Acknowledgements}~: I thank my thesis advisor Olivier Druet for his constant support during the preparation of this paper. I would also like to thank deeply Tristan Rivi\`ere for his valuable comments and remarks on a first draft of the manuscript. 
 
%%%%%%%%%%%%%%%%%%%%%%%%%%%%%%%%%%
\section{Equation of mean curvature in a $3$-Riemannian manifold}
\label{s2}
Here, we compute the equation satisfied by a conformal immersion with respect to its mean curvature. The fact that we consider conformal immersion is very  natural when we look at problems concerning mean curvature. Especially in dimension $2$ where, thanks to the uniformization theorem, on the sphere every metric is conformally equivalent to the standard metric.\\

Let $(\mathcal{N},g)$ be an oriented $3$-Riemanniann manifold, $(\mathcal{M},h)$ be an oriented surface and $f:\mathcal{M} \rightarrow \Sigma\subset\mathcal{N}$ be a {\bf conformal immersion}, that is to say such that $f^{*}(g_{\vert \Sigma}) = e^{2u} h$, where $u\in C^{\infty}(\mathcal{M})$. Hence $\Sigma=f(\mathcal{M})$ is a surface of mean curvature equal to $H$ if
\beq
\label{main1}
 \frac{\partial^{2} f^{j}}{(\partial x^{\alpha})^2}+\Gamma^{j}_{ik}(f) \frac{\partial f^{i}}{\partial x^{\alpha}}\frac{\partial f^{k}}{\partial x_{\alpha}} = 2H(f) \sqrt{\vert g\vert} g^{ij}(f)\nu_{i} \hbox{ for } j\in \{1,2,3\} . \\
\eeq
In arbitrary coordinates, (\ref{main1}) is transformed into 
\beq
\label{main2}
\Delta_{\mathcal{M}} f^{j} - h^{\alpha\beta}\Gamma^{j}_{ik}(f) \frac{\partial f^{i}}{\partial x^{\alpha}}\frac{\partial f^{k}}{\partial x^{\beta}} = -2 H(f) \frac{\sqrt{\vert g\vert} g^{ij}(f)\nu_{i}}{\sqrt{\vert h \vert}} \hbox{ for } j\in \{1,2,3\} , \\
\eeq
where $\Delta_\mathcal{M}$ is the Laplace-Beltrami operator of $(\mathcal{M},h)$. Here we have to notice the fundamental fact that this equation is invariant by a conformal diffeomorphism. That is to say, if $u$ satisfies (\ref{main2}) and $\phi\in \hbox{Conf}(\mathcal{M})$ then $u\circ\phi$ still satisfies (\ref{main2}).\\

{\bf The fundamental example of the Euclidean case : }\\

If we consider a sphere with mean curvature $H$ immersed in $\R^{3}$, we get the so-called equation of H-bubbles :
\beq
\label{ex}
\begin{cases}
\Delta_{\xi} u = -2 H(u)\, u_{x} \wedge u_{y}, \\ 
\langle u_x, u_y\rangle_\xi = 0 \hbox{ and } \Vert u_x \Vert_\xi = \Vert u_y \Vert_\xi.
\end{cases}
\eeq
Here $\xi$ is the standard metric of $\R^3$. This equation, in particular when $H$ is constant, will play a fundamental role in what follows since this is the limit of the general equation when the metric becomes flat. Moreover, thanks to Hopf's theorem, we know that the sphere is the only immersed compact simply connected surface with constant mean curvature in $\R^3$. Hence the sphere provides us a fundamental solution of (\ref{ex}). Hence we need a conformal parametrization of the sphere: this is exactly the purpose of the inverse of the stereographic projection.\\ 

Let $\omega :\R^2 \rightarrow \R^3$ be defined as follows
$$ \omega(x,y)= \frac{1}{1+r^2}\left(\begin{array}{c}2x \\2y \\r^2-1\end{array}\right), $$
where $r^2=x^2+y^2$. This exactely the inverse ogf the stereographic projection with respect to the north pole. Computing the derivatives and their cross product, we get the following useful formulas
$$ \omega_x (x,y)= \frac{2}{(1+r^2)^2}\left(\begin{array}{c}1+(y^2-x^2) \\-2xy \\ 2x\end{array}\right),\ \omega_y(x,y)= \frac{2}{(1+r^2)^2}\left(\begin{array}{c}-2xy\\ 1+(x^2-y^2) \\ 2y\end{array}\right),$$
$$ \omega_x \wedge \omega_y (x,y)= \frac{-4}{(1+r^2)^3}\left(\begin{array}{c}2x \\ 2y\\ 1-r^2\end{array}\right) =  \frac{-4\omega(x,y)}{(1+r^2)^2},$$
\beq
\label{di}
\frac{\vert \nabla \omega\vert^2}{2}= \vert \omega_x \vert^2 = \vert \omega_y \vert^2 =\frac{4}{(1+r^2)^2} \hbox{ and }
\langle \nabla \omega^k ,  \nabla \omega^l \rangle = (\delta_{kl} -   \omega^k   \omega^l) \frac{\vert \nabla \omega \vert^2}{2}.
\eeq
Then we remind a very important result of Brezis and Coron \cite{BC2} which states that the only solution of
$$\Delta u = -2\, u_x \wedge u_y, $$
with bounded energy are exactly, up to a conformal reparametrization, the inverse of the stereographic projection. This result can be seen as a variant of the Hopf's theorem, see \cite{Hopf}, where the hypothesis of conformality is replaced by a bound on the area.   
\begin{lemma}[lemma A.1 of \cite{BC2}]
\label{bcl}
Let $ u \in L^{1}_{loc}(\R^2,\R^3)$ which satisfies 
\beq
\label{eqlim}
\begin{split}
&\Delta u = -2\, u_x \wedge u_y ,\\
& \int_{\R^2} \vert\nabla u \vert^2 dz <+\infty .
\end{split}
\eeq
Then $u$ has precisely the form 
$$u(z) =\omega \left(\frac{P(z)}{Q(z)} \right) + C, $$
where $P$ and $Q$ are polynomial, C is a constant. In addition
$$\int_{\R^2} \vert \nabla u \vert^2\, dz = 8\pi k \hbox{ with }
 k=max\{ deg P, deg\ Q\},$$
 provided that $\frac{P}{Q}$ is irreducible.
 \end{lemma}
 It could be useful to remark that, thanks to (\ref{di}), the gradient of such an $\omega$ satisfies the following formula
\be
\vert \nabla \omega \vert = \frac{2\sqrt{2} \vert P' Q - Q'P\vert}{\vert P\vert^2 +\vert Q\vert^2} . 
\ee
Then we defined a special class of solutions which will be very important in what follows: the spheres which are parametrized only one time.
 \begin{defi} A solution $u$ of (\ref{eqlim}) is said to be simple if
 $$u(z) =\omega \left(\frac{P(z)}{Q(z)} \right) + C, $$
 with  $\frac{P}{Q}$ is irreducible and $max\{ deg P, deg\ Q\}= 1$.
 \end{defi}
In particular, if $u$ is a simple solution of (\ref{eqlim}), then we have 
\beq
\label{oi}
\left\vert \nabla \omega^\eps (x) \right\vert = O \left(\frac{\lambda^\eps}{\vert x- a^\eps\vert^2 +(\lambda^\eps)^2} \right), 
\eeq
where  $u^\eps = u \left( \frac{\, .\, - a^\eps}{\lambda^\eps} \right)$, $a^\eps$ and $\lambda^\eps$ are respectively a sequence of points in $\R^2$ and  a sequence of positive numbers.

%%%%%%%%%%%%%%%%%%%%%%%%%%%%%%%%
 %%%%%%%%%%%%%%%%%%%%%%%%%%%%%%%%%
  
\section{Preliminaries}
\label{s3}
The aim of this section is to remind some basic facts about embedded surfaces in the euclidean space and to use them to give an appropriate formulation of the problem.\\ 

First of all, we give  some classical relations between the diameter, the area and the mean curvature of such embedded surfaces. Then we will give an equivalent of such relations in our Riemannian setting.\\

The following classical  lemma gives a lower bound of the diameter by the inverse of the mean curvature.
\begin{lemma} 
\label{est1}
Let $S$ be a smooth surface of $\R^3$ with mean curvature $H$. Then 
$$2 \leq \delta(S) \sup_{x\in S} \vert H(x) \vert,$$
where $\delta(S)$ is the extrinsic diameter of $M$.
\end{lemma}

{\it Proof of lemma \ref{est1}:}\\

Let $\overline{B(x, r)}$ be the smallest closed ball that enclose $M$. Using a classical maximum principal at $y\in \overline{B(x, r)} \cap S$, we see that $\vert H(y) \vert  \geq \vert H_{S(x,r)}(y) \vert =\frac{1}{r}$, which proves the lemma.\hfill$\square$\\

Then, we remind the Simon's inequality which relates the  diameter to the area and the mean curvature. 
\begin{thm} Let $S$ be a closed connected surface immersed in $\R^3$, then
\beq
\label{Simon}
\delta(S) < \frac{2}{\pi} A(S)^\frac{1}{2} \left( \int_S \vert H \vert^2 d\sigma \right)^\frac{1}{2},
\eeq
where $A$, $H$ and $d\sigma$ are respectively the area, the mean curvature and the volume element of $S$.
\end{thm}
See \cite{Simon} for the original proof and \cite{Topping98} for  the proof with the optimal constant $\frac{2}{\pi}$. Indeed considering a long cylinder ending by spherical caps we see that the constant cannot be improved. Then we obtain, as a by product of (\ref{Simon}), that
\beq
\label{Simon2}
\delta(S) < \frac{2}{\pi} A(S) \sup_{x\in S} \vert H(x) \vert .
\eeq
Such an inequality have been also proved by Bethuel and Rey, see theorem 6.2 of \cite{BethuelRey}. By now, we are in position to give a proof, in Riemannian setting, of the fact that the diameter of an $H$-surface is controlled by the product of the area by the mean curvature. In fact, we need the additional assumption that the diameter is small enough, in order to get a relatively flat geometry. It is false without this additional assumption as shown by a tubular surface around a closed geodesic of  $S^{3}$, see \cite{MazzeoPacard}. 
\begin{lemma} 
\label{estim}
Let $(\mathcal{N},g)$ be a $3$-Riemannian manifold whose injectivity radius admits a positive lower bound and $\Sigma_H \subset \mathcal{N}$ a sequence of connected $H$-surface which satisfies  the following hypothesis
\be
\begin{array}{cc} \begin{cases}
 \delta  (\Sigma_{H}) =o(1) ,\\
 A(\Sigma_{H}) =O\left(\frac{1}{H^2}\right).
\end{cases} & \hbox{ as } H\rightarrow +\infty,  \end{array}
\ee 
Then we get the following estimate 
$$ \frac{1}{K H}\leq  \delta(\Sigma_H) \leq \frac{K}{H} ,$$
where $K$ is a positive constant.
\end{lemma}

{\it Proof of lemma \ref{estim} :}\\

Let $c_{H}\in \Sigma_{H}$, for $H$ large enough, we can assume that $\Sigma_H\subset B(c_H,\delta)$ where $\delta$ is smaller than the injectivity radius of $\mathcal{N}$. Then we rescale  the exponential chart centered in $c_H$ by a factor $\frac{1}{\delta_H}$, where $\delta_H = \delta(\Sigma_H)$. We get a sequence of surface~$\tilde{\Sigma}_{H}$ of $(\R^3, g_H)$ with diameter $1$, here $g_H$ is the rescale metric. The mean  curvature of $\tilde{\Sigma}_H$, computed with respect to $g_H$, is equal to $\tilde{H}= \delta_H H$. Remarking that  $g_H$ converges uniformly to $\xi$ on every compact and that $\tilde{\Sigma}_H\subset \overline{B(0,2)}$ we get, thanks to lemma \ref{est1}, for $H$ large enough that
$$ \tilde{H}= \delta_H H\geq 1,$$
which proves the left hand-side inequality.\\
In the other hand, thanks to (\ref{Simon2}), for $H$ large enough, we have
$$ 
\tilde{H}=\delta_H H \geq \frac{(\delta_H H)^2}{C}  ,
$$  
where $C$ is the positive constant. This achieves the proof of the lemma.\hfill$\square$\\

Hence, we immediately see that our assumption (\ref{H}) is equivalent to assuming that $\Sigma_H$ satisfies
\beq
\tag{H'}
\label{H'}
\begin{array}{cc} \begin{cases}
\frac{1}{CH}\leq \delta(\Sigma_{H}) \leq \frac{C}{H}\\
 A(\Sigma_{H}) \leq \frac{C}{H^2}
\end{cases} & \hbox{ as } H\rightarrow +\infty,  \end{array}
\eeq 
where $C$ is a positive constant.\\

In order to look more precisely at our $H$-surfaces we need some coordinates. In particular we have to choose a center of chart. For that purpose we fix an arbitrary point $c_H$ of $\Sigma_H$ as a centre of chart. Up to a subsequence, $\Sigma_{H}\rightarrow p_{\infty}$ as $H\rightarrow +\infty$. Of course $c_H\rightarrow p_\infty$ as $H\rightarrow +\infty$. From now on, we look at $\Sigma_H$ in the exponential chart centered at $c_H$. Then we rescale this chart by a factor $H$ with respect to $0$ and we replace the variable $H$ by $\frac{1}{\eps}$. Hence we get a new sequence of immersed spheres $(\Sigma_{\eps}) \subset (\R^3,g_\eps)$ with constant mean curvature $1$, where $g_{\eps}$ is the rescaled metric: $g_\eps(y)(u,v)= g(\eps y) (\eps u,\eps v)$. Moreover, $\Sigma_{\eps}$ satisfies the following additional assumption  
\beq
\tag{$H''$}
\label{H''}
\begin{cases}
A(\Sigma_{\eps}) \leq C ,\\
\Sigma_{\eps} \subset B(0, C),
\end{cases}
\eeq
where $C$ is a positive constant.\\

Finally, let  $u^\eps$  be a parametrization of $\Sigma_\eps$ from $(S^2,h)$ to $(\R^3,g_\eps)$. Up to a diffeomorphism of the sphere we can assume this parametrization to be conformal. Indeed $(u^\eps)^*({g_\eps}_{\vert \Sigma_\eps})$ is in the conformal class of the standard metric, since there is only one conformal class on $S^2$. Hence, let $\phi^\eps \in \hbox{Diff}(S^2)$ such that $(\phi^\eps)^*((u^\eps)^*({g_\eps}_{\vert \Sigma_\eps}))$ is pointwise conformal to $h$, then $u^\eps \circ \phi^\eps$ is our conformal parametrization  from $(S^2,h)$ to $(\R^3,g_\eps)$. Up to replace $u^\eps$ by  $u^\eps \circ \phi^\eps$, $u^\eps$ satisfies, in any conformal coordinates, the following equations
\beq
\label{main}
\begin{cases}
 \Delta_{S^2} u^{\eps} - (\Gamma^{j}_{ik})_\eps(u^\eps) \langle \nabla (u^i)^{\eps}, \nabla (u^k)^{\eps}\rangle_h = -2 \frac{\sqrt{\vert g_\eps \vert} g_\eps^{ij} ((u^\eps)_x \wedge (u^\eps)_y)_i}{\sqrt{\vert h \vert}} \\
  \Vert u^\eps\Vert_\infty \leq C\\
  \Vert \nabla u^\eps\Vert_2 \leq C ,
 \end{cases} 
\eeq
where $(\Gamma^{j}_{ik})_\eps$ are the Christoffel symbols of $g_\eps$ and $C$ is a positive constant. This equation is totally invariant under any conformal diffeomorphism of the sphere. But as we will remind it, the group of conformal diffeomorphism of the $S^2$, $\mathrm{Conf}(S^2)$, is not compact. Hence it could be interesting to fix our parametrization once and for all. Of course, there is no canonical choice. We choose to rescale the highest bubble around the north pole and so to send the remainder around the south pole.\\

 But before making this rescaling, let us defined a dilatation on the sphere $S^2$. For any $Q\in S^2$, let  $\pi_{Q} : S^2 : \rightarrow \C$ the associate  stereographic projection (here $(\R Q)^\bot$ is identified to  $\C$) and for any $t\in [1,+\infty [$. Let $\tau_{t} : \C \rightarrow \C $ as $\tau_{t}(z)=tz$, we set
$$ \Phi_{Q,t}=\pi^{-1}_{Q} \circ \tau_{t} \circ \pi_{Q} .$$
Hence,  we are in position to fix the parametrization of 
$u^\eps$. Let $a^\eps$  and $\lambda^\eps$ be such that 
$$\vert \nabla u^\eps (a^\eps)\vert  = \frac{1}{\lambda^\eps} = \sup_{S^2} \vert\nabla u^\eps\vert .$$
Up to compose $u^\eps$ with a rotation of $S^2$, we can also  assume that $a^\eps=N$. Then we replace $u^\eps$ by  $u^\eps \circ \Phi_{N, \lambda^\eps}$ and we easily check that $\nabla u^\eps$ is bounded on every compact subset of $S^2\setminus\{ S\}$. Moreover, thanks to the conformal invariance of our problem, $u^\eps $ still satisfies (\ref{main}). Hence, thanks to standard elliptic theory, see \cite{GT}, there exist a subsequence of $u^\eps$ (still denoted  $u^\eps)$ and $u^0\in C^2(S^2 \setminus\{ S\})$  such that 
\beq
u^\eps \rightarrow u^0 \hbox{ in } C^2_{loc}(S^2 \setminus\{ S\}) ,
\eeq
If we set  $\omega^0 = u^0 \circ \pi_N^{-1}$, then $\omega^0 \in C^2 (\R^2\setminus \{0\})$ and  satisfies 
$$  \Delta \omega^0 = -2\,  \omega^0_x \wedge \omega^0_y \hbox{ on } \R^2\setminus \{0\}.$$
Then, thanks to the conformal invariance  of $\Vert\nabla\, .\, \Vert_2$, we have
$$\Vert \nabla \omega^0 \Vert_2 \leq \liminf_{\eps\rightarrow 0} \Vert \nabla u^\eps \circ \pi_N^{-1} \Vert_2= \liminf_{\eps\rightarrow 0} \Vert \nabla u^\eps\Vert_2<+\infty .$$
Hence $\omega^0$ is a solution of (\ref{eqlim}) and $\omega^{0}$ is non trivial since $\vert \nabla u^0 (N)\vert=1$. Moreover $\vert \nabla \omega^0\vert$ has a maximum in $\R^2$, let $a_0\in\R^2$ be a point where  $\vert \nabla \omega^0\vert$ achieves its maximum.\\ 

Finally, up to replace $u^\eps$ by $u^\eps \circ \pi_N^{-1}$,  $u^\eps$ satisfies
\beq
\label{main3}
\begin{cases}
 \Delta_{\xi} u^{\eps} - (\Gamma^{j}_{ik})_\eps(u^\eps) \langle \nabla (u^i)^{\eps}, \nabla (u^k)^{\eps}\rangle_\xi = -2 \sqrt{\vert g_\eps\vert} g_\eps^{ij} ((u^\eps)_x \wedge (u^\eps)_y)_i\\
 u^\eps_0 \rightarrow \omega^0 \hbox{ in } C^2_{loc} (\R^2\setminus\{ -a_0\}) \\
  \Vert u^\eps\Vert_\infty \leq C\\
  \Vert \nabla u^\eps\Vert_2 \leq C .
 \end{cases} 
\eeq
where  $u^\eps_0=u^\eps(z+a_0)$, $a_0\in\R^2$ and $\omega^0$ is a non trivial solution of (\ref{eqlim}) such that  $\vert \nabla \omega^0\vert$ achieves its maximum at $0$.\\

Now on and until the end,  $u^\eps$ is seen as a map from $\R^2$ to $\R^3$.
%%%%%%%%%%%%%%%%%%%%%%%%%%%%%%%%%%%%%%%%%%%%%%%%%%%%%%%% 

\section{Decomposition of $u^\eps$ as sum of bubbles.}
\label{s4}
The aim of this section consists in two steps. First we will show that $\Sigma_\eps$  converges to a sum of round spheres. Then we will adjust these round spheres to the geometry of our manifold. All of this will be sum up at the end of this section.\\

Such a decomposition has already been observed by  Brezis and Coron in \cite{BC2} where they notably give an $H^1$-decomposition for approached solution of the mean curvature equation on the disk. Here we give a result in the same spirit, replacing the $H^1$ by $C^2_{loc}$. The method used  have been intensively used for the Yamabe equation and then generalized to critical elliptic systems, see \cite{DruetHebeyRobert}, \cite{DruetHebey05} and \cite{DruetHebey09}.
\begin{thm}
\label{de}
Let $u^\eps$ be a sequence of $C^2$-solutions of (\ref{main3}). Then, there exist $p \in \N$ and 
\begin{enumerate}[(i)]
\item   $\omega^1,\dots,\omega^p $ simple solutions of (\ref{eqlim}) such that $\vert \nabla \omega^i \vert$ has a maximum at $0$, 
\item  $ a^\eps_1, \dots, a^\eps_p$ sequences of $\R^2$ which all converge to $0$, and
\item  $\lambda^\eps_1 , \dots , \lambda^\eps_p$  sequences of positive numbers such that $\ds \lim_{\eps \rightarrow 0} \lambda^\eps_i =0$, 
\end{enumerate} 
such that, for a subsequence of $u^\eps$ (still denoted $u^\eps$ ) the following assertions hold 
\beq
\label{A}
\tag{A}
 u^\eps_i \rightarrow \omega^i \hbox{ in } C^2_{loc}(\R^2\setminus S_i) \hbox{ as } \eps \rightarrow 0 \hbox{ for all } 1\leq i\leq p,
 \eeq
 where $u^\eps_i= u^\eps (\lambda^\eps_i\ . + a_i^\eps )$ and $\ds S_i= \lim_{\eps \rightarrow 0} \left\{ \frac{a_j^\eps-a_i^\eps}{\lambda_i^\eps}\hbox{ s.t. } j\in \{1,\dots,p\}\setminus \{ i\}\right\} $.

\beq
\label{ortho}
\tag{B}
 \frac{d_{i}^\eps(a_j^\eps)}{\lambda^\eps_j} + \frac{d_{j}^\eps(a_i^\eps)}{\lambda^\eps_i} \rightarrow +\infty  \hbox{ for all }i \not= j ,
\eeq
 where $d_{i}^\eps(x)=\sqrt{(\lambda^\eps_i)^2 +\vert a_i^\eps -x\vert^2}$ for $1\leq i \leq p$ and $d_0(x)=\sqrt{1+\vert x -a_0\vert^2}$.\\
 
With the additional properties that 
\beq
\label{faible}
\tag{C}
 \lim_{\eps \rightarrow 0} \sup_{x\in \R^2}  \left(\min_{0 \leq i \leq p} d_i^\eps (x)\right) \left\vert \nabla  \left(u^\eps-\sum_{i=0}^p  \omega^\eps_i\right)  (x) \right\vert =0
\eeq
and 
\beq
\tag{D}
\label{h1dec}
\left\Vert\nabla \left( u^\eps  -\sum_{i=0}^p \omega^\eps_i \right) \right\Vert_2 \rightarrow 0  \hbox{ as } \eps\rightarrow 0,
\eeq 
where $\omega^\eps_i = \omega_i\left(\frac{\ .\ -a_i^\eps}{\lambda_i^\eps} \right)$ and $(a_0^\eps,\lambda^\eps_0)=(a_0,1)$.
\end{thm}
 When  there is just one bubble, that is to say when  $p=0$, the conclusion limits to
\be
u^\eps \rightarrow \omega^0 \hbox{ in } C^2(\R^2) \hbox{ as } \eps \rightarrow 0.
\ee
 
{\it Proof of theorem \ref{de} :}\\

We are going to extract the bubbles by induction and the process will stop thanks to our uniform bound on the area of $\Sigma_\eps$.\\

{\bf For $k\geq 0$ let $(P_k)$ be  the following assertion :}\\

There exist
\begin{enumerate}[(i)]
\item  $\omega^0,\dots,\omega^k $ non trivial solutions of (\ref{eqlim}) such that $\vert \nabla \omega^i \vert$ has its maximum at $0$, 
\item  $ a^\eps_0, \dots, a^\eps_k$  bounded sequences of $\R^2$ such that $\ds \lim_{\eps\rightarrow 0} a_i^\eps =0$ for $1\leq i \leq k$
, and
\item  $\lambda^\eps_0, \dots,\lambda^\eps_k$ bounded sequences of positive numbers such that $\ds \lim_{\eps \rightarrow 0} \lambda^\eps_i =0$ for $1\leq i \leq k$, 
\end{enumerate} 
such that, for a subsequence of $u^\eps$ (still denoted $u^\eps$) the following assertions hold

\beq
\tag{$A_k$}
\label{ui}
 u^\eps_i \rightarrow \omega_i \hbox{ in } C^2_{loc}(\R^2\setminus S_i) \hbox{ as } \eps \rightarrow 0 \hbox{ for all },
 \eeq
 where $u^\eps_i= u^\eps (\lambda^\eps_i\ .\  -a_i^\eps)$ and $\ds S_i= \lim_{\eps \rightarrow 0} \left\{ \frac{a_j^\eps -a_i^\eps}{\lambda_i^\eps}\hbox{ s.t. } j\in \{0,\dots,k\}\setminus \{ i\}\right\} .$
 \beq
 \tag{$B_k$}
\label{ortho1}
 \frac{d_{i}^\eps(a_j^\eps)}{\lambda^\eps_j} + \frac{d_{j}^\eps(a_i^\eps)}{\lambda^\eps_i} \rightarrow +\infty \ \forall i \not= j , \hbox{ as } \eps \rightarrow 0,
\eeq
 where $d_{i}^\eps(x)=\sqrt{(\lambda^\eps_i)^2 +\vert a_i^\eps -x\vert^2}$.\\

{\bf Claim 1: if $(P_k)$ holds  for some $k\geq 0$ then either $(P_{k+1})$ holds or
\beq
\lim_{\eps \rightarrow 0} \sup_{x\in \R^2}  \left(\min_{1\leq i \leq k} d_i^\eps (x)\right) \left\vert \nabla \left(u^\eps -\sum_{i=0}^k  \omega^\eps_i\right) (x)\right\vert =0,
\eeq
where $\omega^\eps_i = \omega_i\left(\frac{\, . \, -a_i^\eps}{\lambda_i^\eps} \right)$.
}\\

{\it Proof of Claim 1 :}\\

In order to prove this claim, we assume that $(P_k)$ holds and that there exists $\gamma_0>0$ such that
\beq
\label{gam0}
 \sup_{z\in \R^2}  \left(\min_{0\leq i \leq k} d_i^\eps (z)\right) \left\vert \nabla \left(u^\eps-\sum_{i=0}^k  \omega^\eps_i\right) (z) \right\vert \geq \gamma_0 \hbox{ for all } \eps>0.
\eeq
We need to prove that $(P_{k+1})$ holds. Let $a^\eps_{k+1} \in \R^2$ be such that
\be
\left(\min_{0\leq i\leq k} d_i^\eps (a_{k+1}^\eps) \right) \left\vert \nabla  \left(u^\eps -\sum_{i=0}^k  \omega^\eps_i\right)  (a_{k+1}^\eps) \right\vert= \sup_{z\in \R^2} \left( \min_{0\leq i \leq k} d_i^\eps (z) \right) \left\vert \nabla  \left(u^\eps- \sum_{i=0}^k  \omega^\eps_i\right)  (z) \right\vert .
\ee 
The fact that the supremum is achieved is a consequence of our assumptions. Indeed, thanks (\ref{main3}), we get
$$\left\vert \nabla \left(u^\eps-\sum_{i=0}^k  \omega^\eps_i\right) (z) \right\vert =O\left( \frac{1}{1+\vert z\vert^2}\right) \hbox{ as } z \rightarrow +\infty ,$$ 
which proves that the maximum is achieved. Now we define $\lambda^\eps_{k+1}$ by the equation
\be
\left\vert \nabla\left(u^\eps -\omega^0 -\sum_{i=1}^k  \omega^\eps_i\right) (a_{k+1}^\eps)\right\vert  = \frac{1}{\lambda_{k+1}^\eps} .
\ee
Always thanks to (\ref{main3}) and the assumptions about the  $a_i^\eps$ and the  $\lambda_i^\eps$ , we remark that
\be
\left\vert \nabla\left(u^\eps -\sum_{i=0}^k  \omega^\eps_i\right)\right\vert \rightarrow 0  \hbox{ on } \R^2\setminus \{ 0 \} ,
\ee
then $(a^\eps_{k+1})$ converges to $0 $. Then, if $k\geq 1$ we have 
\be
\min_{0\leq i\leq k} d_i^\eps (a_{k+1}^\eps) \rightarrow 0 \hbox { as } \eps\rightarrow 0 ,
\ee
and
\beq
\label{l} 
\lambda_{k+1}^\eps \rightarrow 0 \hbox { as } \eps\rightarrow 0 .
\eeq
In fact, (\ref{l}) is also true when $k=0$. Indeed, else $u^\eps -\omega_0^\eps$ would be uniformly bounded in $C^1(\R^2)$ and hence converge to $0$ on the whole plane which contradicts (\ref{gam0}).

Now there are two cases to consider.\\

{\bf First case :
\beq
\label{c1}
 \lim_{\eps\rightarrow 0} \frac{\displaystyle \min_{0\leq i \leq k} d_i^\eps (a_{k+1}^\eps)}{\lambda_{k+1}^\eps}  = +\infty .
\eeq
}

In this case, ($B_{k+1}$) is automatically satisfied. Now, we set $u^\eps_{k+1} = u^\eps( \lambda_{k+1}^\eps \, . \,  +a_{k+1})$. Let $z\in  \R^2$, we get that 
\beq
\label{c2}
\begin{split}
\vert \nabla u^\eps_{k+1}(z)\vert &=  \lambda_{k+1}^\eps \vert \nabla u^\eps ( \lambda_{k+1}^\eps z + a_{k+1}^\eps)\vert \\
&\leq  \lambda_{k+1}^\eps \left\vert \nabla \left(u^\eps - \sum_{i=0}^k  \omega^\eps_i\right)(  \lambda_{k+1}^\eps z + a_{k+1}^\eps)\right\vert \\
&+  \lambda_{k+1}^\eps \left\vert \nabla \left(\sum_{i=0}^k  \omega^\eps_i\right)( \lambda_{k+1}^\eps z + a_{k+1}^\eps)\right\vert .
\end{split}
\eeq
Thanks to (\ref{oi}) and (\ref{c1}), we easily see that 
\beq
\label{c3}
\begin{split}
&\lambda_{k+1}^\eps \left\vert \nabla \left(\sum_{i=0}^k  \omega^\eps_i\right)(\lambda_{k+1}^\eps z +a_{k+1}^\eps )\right\vert =o(1) ,\\
&\hbox{ and }\\
& \lim_{\eps\rightarrow 0}\lambda_{k+1}^\eps \left\vert \nabla u^\eps ( a_{k+1}^{\eps}) \right\vert =1 .
\end{split} 
\eeq
Then using the definition of $a_{k+1}^\eps$, (\ref{c1}), (\ref{c2}) and (\ref{c3}) we have
\beq
\label{c4}
\vert \nabla u^\eps_{k+1}(z)\vert  \leq \frac{\min_{i} d_i^\eps (a_{k+1}^\eps)}{\min_{i} d_i^\eps (\lambda_{k+1}^\eps z+ a_{k+1}^\eps)} +o(1)=1+o(1) .
 \eeq
Then  $\vert \nabla u^\eps_{k+1}\vert$ is bounded on every compact subset of $\R^2$. Moreover thanks to the conformal invariance of our problem, $u^\eps_{k+1}$ still satisfies (\ref{main3}). Hence, thanks to standard elliptic theory, see \cite{GT}, there exist a subsequence of $u^\eps$ (still denotes  $u^\eps)$ and $\omega^{k+1}\in C^2(\R^2)$  such that 
\be
u^\eps_{k+1} \rightarrow \omega^{k+1} \hbox{ in } C^2_{loc} (\R^2) 
\ee
and 
\be
\Delta \omega^{k+1} = -2\,  \omega^{k+1}_x \wedge \omega^{k+1}_y \hbox{ on } \R^2. 
\ee
Moreover, thanks to the conformal invariance  of $\Vert\nabla .\, \Vert_2$, we have
$$\Vert \nabla \omega^{k+1} \Vert_2 \leq \liminf_{\eps\rightarrow 0} \Vert \nabla u^\eps_{k+1}\Vert_2= \liminf_{\eps\rightarrow 0}\Vert \nabla u^\eps\Vert_2<+\infty .$$
Then, thanks to lemma \ref{bcl}, $\omega^{k+1}$ is a solution of (\ref{eqlim}) on $ \R^2$  and $\omega^{k+1}$ is non-trivial since $\vert \nabla \omega^{k+1}(0)\vert=1$. Finally, thanks to (\ref{c3}) and (\ref{c4}), we easily see that $\vert \nabla \omega^{k+1}\vert$ has a maximum at $0$. This achieves the proof of the fact that $(P_{k+1})$ holds in the first case.\\

{\bf Second case :
\beq
\label{29}
 \lim_{\eps\rightarrow 0} \frac{\ds \min_{0\leq i \leq k} d_i^\eps (a_{k+1}^\eps)}{\lambda_{k+1}^\eps}  = \gamma >0 .
\eeq
}
In that case we necessary get $k>0$. \\

First of all, we need to prove that  ($B_{k+1}$) holds. If it doesn't hold, up to a subsequence, there exists $1\leq i_0 \leq k$ such that
\beq
\label{c12} 
d_{k+1}^\eps (a_{i_0}^\eps)=O(\lambda_{i_0}^\eps) \hbox{ and } d_{i_0}^\eps (a_{k+1}^\eps) =O(  \lambda_{k+1}^\eps) . 
\eeq
From the one hand,  (\ref{c12}) gives that
\beq
\label{c14}
\begin{split}
\frac{\lambda_{k+1}^\eps}{\lambda_{i_0}^\eps} \rightarrow c \hbox{ as } \eps\rightarrow 0
\hbox{ and }
\vert a_{i_0}^\eps - a_{k+1}^\eps \vert = O(\lambda_{i_0}^\eps),
\end{split}
\eeq
where $c$ is a positive constant. From the other hand, thanks to (\ref{ui}) and (\ref{ortho1}), we have 
\beq
\label{c13}
 \nabla \left(\left(u^\eps-\sum_{i=0}^k  \omega^\eps_i\right) (\lambda_{i_0}^\eps\ . + a_{i_0}^\eps)\right) \rightarrow 0 \hbox{ in }C^1_{loc} (\R^2\setminus\{S_{i_0}\}). 
 \eeq
Then, thanks to (\ref{29}) and (\ref{c14}), we necessary get  that 
$$ d\left(\frac{a_{k+1}^\eps - a_{i_0}^\eps}{\lambda_{i_0}^\eps}, S_{i_0}\right)= o(1).$$ 
Then there exists $j\in \{0,\dots,k\}\setminus\{i_0\}$ such that 
$$\left\vert \frac{a_{k+1}^\eps - a_{j}^\eps}{\lambda_{i_0}^\eps}\right\vert= o(1).$$ 
Then, thanks to (\ref{29}) and (\ref{c14}), for $\eps$ small enough,  we get 
$$ \frac{\lambda_j^\eps}{\lambda_{k+1}^\eps} \geq \frac{\gamma}{2} ,$$
and, thanks to (\ref{c14}), for $\eps$ small enough, we get
$$
\frac{\lambda_j^\eps}{\lambda_{i_0}^\eps}  \geq \frac{\gamma}{4c} ,
$$
But, since $\frac{a_{i_0}^\eps -a_j^\eps}{\lambda_{i_0}^\eps} =O(1)$ and that $i_0$ and $j$ satisfies (\ref{ortho1}), we have
$$
\lambda_{i_0}^\eps =o(\lambda_j^\eps) .
$$
Hence for every $j$ such that  $\frac{a_{k+1}^\eps -a_j^\eps}{\lambda_{i_0}^\eps} =o(1)$  we have
$$
\lambda_{i_0}^\eps =o(\lambda_j^\eps) .
$$
In particular, thanks to (\ref{oi}), there exists $\delta >0$ such  that for every $z\in B(0,\delta)$ we get that
$$\lambda_{i_0}^\eps \vert \nabla \omega_i^\eps (a_{k+1}^\eps + z \lambda_{i_0}^\eps)\vert =o(1) \hbox{ for every } i\not= i_0$$
 Then we easily get that 
$$\lambda_{i_0}^\eps \vert \nabla u^\eps \vert = O(1) \hbox{ on } B(a_{k+1}^\eps, \delta \lambda_{i_0}^\eps) .$$
Hence thanks to elliptic theory, up to a subsequence, we see that 
$$ \vert \nabla (u^\eps_{i_0} -\omega^{i_0})(a^\eps_{k+1})\vert \rightarrow 0,$$
which leads to
$$
\lambda_{i_0}^\eps \left\vert\nabla \left(\left(u^\eps-\sum_{i=0}^k  \omega^\eps_i\right) (a_{k+1}^\eps)\right)\right\vert \rightarrow 0, 
$$
which, thanks to (\ref{c14}), is a contradiction with (\ref{29}) and proves ($B_{k+1}$).\\ 

Now, we set $u^\eps_{k+1} = u^\eps(\lambda_{k+1}^\eps\ .+ a_{k+1}^\eps)$. Let $z\in  \R^2\setminus\{S_{k+1}\}$, we get that 
\beq
\label{c22}
\begin{split}
\vert \nabla u^\eps_{k+1}(z)\vert &=  \lambda_{k+1}^\eps \vert \nabla u^\eps (\lambda_{k+1}^\eps z+ a_{k+1}^\eps)\vert\\
&\leq  \lambda_{k+1}^\eps \left\vert \nabla \left(u^\eps - \sum_{i=0}^k  \omega^\eps_i\right)(\lambda_{k+1}^\eps z+ a_{k+1}^\eps)\right\vert \\
&+  \lambda_{k+1}^\eps \left\vert \nabla \left(\sum_{i=0}^k  \omega^\eps_i\right)(\lambda_{k+1}^\eps z + a_{k+1}^\eps)\right\vert .
\end{split}
\eeq
Thanks to (\ref{oi}) and  (\ref{29}), we easily see that 
\beq
\label{c32}
\lambda_{k+1}^\eps \left\vert \nabla \left(\sum_{i=0}^k  \omega^\eps_i\right)(\lambda_{k+1}^\eps\ .+ a_{k+1}^\eps)\right\vert =O\left(\frac{1}{d(z,S_{k+1})}\right) .
\eeq
Then using the definition of $a_{k+1}^\eps$, (\ref{c22}) and (\ref{c32}) we have
\beq
\label{c33}
\vert \nabla u^\eps_{k+1}(z)\vert  \leq \frac{\min_{i} d_i^\eps (a_{k+1}^\eps)}{\min_{i} d_i^\eps (\lambda_{k+1}^\eps z + a_{k+1}^\eps)} +  O(\frac{1}{d(z,S_{k+1})})=O(\frac{1}{d(z,S_{k+1})}) .
 \eeq
Then  $\vert \nabla u^\eps_{k+1}\vert$ is bounded on every compact subset of $\R^2\setminus\{S_{k+1}\}$. Moreover thanks to the conformal invariance of our problem, $u^\eps_{k+1}$ still satisfies (\ref{main3}). Hence, thanks to standard elliptic theory, see \cite{GT}, there exists a subsequence of $u^\eps$ (still denotes  $u^\eps)$ and $\omega^{k+1}\in C^2(\R^2\setminus S_{k+1})$  such that 
\be
u^\eps_{k+1} \rightarrow \omega_{k+1} \hbox{ in } C^1_{loc} (\R^2\setminus S_{k+1})
\ee
and 
$$  \Delta \omega_{k+1} = -2  \omega^{k+1}_x \wedge \omega^{k+1}_y \hbox{ on } \R^2\setminus S_{k+1}.  $$
Moreover, thanks to the conformal invariance  of $\Vert\nabla\, . \Vert_2$, up to extraction, we have
$$u^\eps_{k+1}\rightharpoonup \omega^{k+1} \hbox{ in } L^2 (\R^2) $$
and 
$$\Vert \nabla \omega^{k+1} \Vert_2 \leq \liminf_{\eps\rightarrow 0} \Vert \nabla u^\eps_{k+1}\Vert_2= \liminf_{\eps\rightarrow 0} \Vert \nabla u^\eps\Vert_2<+\infty .$$
Then, thanks to lemma \ref{bcl}, $\omega^{k+1}$ is a solution of (\ref{eqlim}) on $ \R^2$. Then we want to show that $\omega^{k+1}$ is non-trivial. This is obvious if $0\not\in S_{k+1}$, since in this case we get  $\vert \nabla \omega^{k+1}(0)\vert=1$. But for every $i_0$ such that 
$$\frac{\vert a_{i_0}^\eps - a_{k+1}^\eps\vert}{\lambda_{k+1}^\eps} =o(1),$$
 thanks to (\ref{29}) and ($B_{k+1}$), we get
$$ \lambda_{i_0}^\eps =o(\lambda_{k+1}^\eps) .$$
Then mimiking the argument of the proof of ($B_{k+1}$) we prove that 
$$ \nabla u^\eps_{k+1} \rightarrow \nabla \omega^{k+1} \hbox{ on } B(0, \delta),$$
where $\delta>0$. Which leads to $\vert \nabla \omega^{k+1}(0)\vert=1$ and proves that $\omega^{k+1}$ is non-trivial.\\ 

Finally $\vert \nabla \omega^{k+1}\vert$ achieves his maximum at $a_{k+1}\in \R^2$, then up to replace $a^\eps_{k+1}$ by $a^\eps_{k+1} + \lambda_{k+1}^\eps a_{k+1}$, the conclusion still holds with a new $\omega^{k+1}$ such that  $\vert \nabla \omega^{k+1}\vert$ achieves his maximum at $0$. This proves ($P_{k+1}$) in the second case. The study of these two cases ends the proof of claim 1.\\

Then, before proving the theorem, we need to prove a claim about the growth of the energy of such a decomposition.\\
  
{\bf Claim 2: Let $k\in \N$ and  
\begin{enumerate}[(i)]
\item   $\omega^0,\dots,\omega^k $ non trivial solution  of  (\ref{eqlim}), 
\item   $ a^\eps_0, \dots, a^\eps_k $ bounded sequences $\R^2$ , and
\item  $\lambda^\eps_0, \dots,\lambda^\eps_k,$ bounded sequences of positive numbers, 
\end{enumerate} 
such that,  with $u^\eps$,  they satisfy ($P_k$). Then
$$ \liminf_{\eps \rightarrow 0} \Vert \nabla u^\eps \Vert_2^2 \geq \sum_{i=0}^k \Vert \nabla \omega^i \Vert^2_2 \geq 8\pi (k+1) .$$
}\\
{\it Proof of claim 2 :}\\

Indeed let $R$ be a real positive number, then, thanks to (\ref{ortho1}), for $\eps$ small enough, we get 
$$
\int_{\R^2} \vert \nabla u^\eps \vert^2   dz \geq \sum_{i=0}^k
\int_{B(a^\eps_i, R\lambda_i^\eps)\setminus \Omega_i^\eps (R)} \vert \nabla u^\eps\vert^2 dz,
$$
where $\ds \Omega_i^\eps (R) = \cup_{j\not= i} B(a^\eps_j, R \lambda_j  ^\eps)$. Then, thanks to ($A_k$), we get 
\beq
\begin{split}
\int_{\R^2} \vert \nabla u^\eps \vert^2   dz
 &\geq   \sum_{i=0}^k \int_{B(0, R)\setminus \Omega_i (R)} \vert \nabla \omega_i \vert^2 dz + \delta_{\eps,R} \\
 & \geq 8\pi (k+1) + \delta_{\eps,R}
\end{split}
\eeq
where $\ds \Omega_i (R)= \cup_{x\in S_i} B(x,\frac{1}{R})$ and $\ds \lim_{R\rightarrow +\infty } \lim_{\eps \rightarrow 0} \delta_{\eps,R}=0$.\hfill$\square$\\

{\bf Proof of the theorem :}\\

Since $u^\eps$ satisfies (\ref{main3}), we see that $(P_0)$ holds and we set $\lambda_0^\eps=1$. Then we can start our extraction. Indeed, thanks to claim 1 and 2 and the fact that $\Vert\nabla u^\eps \Vert_2$ is finite,  there exists $k\in \N$ such that ($P_k$) is satisfied and
\beq
\label{ii}
\lim_{\eps \rightarrow 0} \sup_{x\in \R^2}  \left(\min_{0\leq i\leq k} d_i^\eps (x)\right) \left\vert \nabla \left(u^\eps -\sum_{i=0}^k  \omega^\eps_i\right) (x)\right\vert =0,
\eeq
where $\omega^\eps_i = \omega_i\left(\frac{\, . + a_i^\eps}{\lambda_i^\eps} \right)$. Which proves that (\ref{A}), (\ref{ortho}) and (\ref{faible}) holds. It remains to  prove (\ref{h1dec}). Let 
\be
R^\eps =u^\eps -\sum_{i=0}^k \omega_i^\eps
\ee 
and let us assume for contradiction that there exists $\delta>0$ such that 
\be
\Vert \nabla R^\eps \Vert_2 \geq \delta .
\ee
Then we are going to extract a new bubble and prove this contradicts (\ref{ii}). Here we follow the method developed in \cite{BC2}.\\

First we introduce the concentration function
\be
C^\eps (t) = \sup_{x\in \R^2} \int_{B(x,t)} \vert \nabla R^\eps \vert^2 dz.
\ee
In fact this supremum is a maximum, since  $R^\eps$ is  in $L^2(\R^2)$. Moreover, each $C^\eps$ is continuous, increasing in $t$, $C^\eps(0)=0$ and, thanks to ($P_k$), $C^\eps(1)\geq \frac{C^\eps(\infty)}{2}\geq \frac{\delta}{2}$, for $\eps$ small enough. We fix $\nu$ such that 
\be
0 < \nu < \min \left\{\frac{1}{2C_0}, \frac{\delta}{2} \right\} ,
\ee
where $C_0$ is the constant involved in lemma \ref{wente3}. Hence there exists $a^\eps \in \R^2$ and $\lambda^\eps >0$ such that   
\be
C^\eps (\lambda^\eps) =  \int_{B(a^\eps,\lambda^\eps)} \vert \nabla R^\eps \vert^2 dz= \nu.
\ee
Of course, thanks to ($P_k$), we know that 
\be
a^\eps \rightarrow 0 \hbox{ and } \lambda^\eps \rightarrow 0 \hbox{, as } \eps \rightarrow 0. 
\ee
Then we rescale at $a^\eps$, setting $\tilde{f} =f(\lambda^\eps\, .\, +a^\eps)$, and we get 
\be
\int_{\R^2} \vert \nabla \tilde{R}^\eps \vert^2 dz = \Vert \nabla R^\eps \Vert_2^2 \leq C,
\ee
and
$$\Vert \nabla \tilde{R}^\eps \Vert_\infty \leq C, $$
where $C$ is a positive constant. Moreover, thanks to (\ref{main3}), $\tilde{R}^\eps$ satisfies
\be
\begin{split}
\Delta \tilde{R}^\eps &= -2\, \tilde{R}^\eps_x \wedge \tilde{R}^\eps_y +O \left(\sum_{i=0}^k \vert \nabla \tilde{\omega}_i^\eps\vert \left( \sum_{ j\not= i}  \vert \nabla \tilde{\omega}_j^\eps\vert  + \vert \nabla \tilde{R}^\eps\vert\right)  \right) \\
&+ O(\eps^2 \vert \nabla\tilde{u}^\eps \vert^2 ) .
\end{split}
\ee
From the other hand, we get, thanks to (\ref{ortho1}), that 
\be
\vert \nabla \tilde{\omega}_i^\eps\vert \vert \nabla \tilde{\omega}_j^\eps\vert \rightarrow 0 \hbox{ in } L^1_{loc}(\R^2)
\ee
and, thanks to (\ref{ii}), we get that
\be
\vert \nabla \tilde{\omega}_i^\eps\vert \vert \nabla \tilde{R}^\eps\vert  \rightarrow 0 \hbox{ in } L^1_{loc}(\R^2) .
\ee
Hence we get that 
\be
 \Delta \tilde{R}^\eps= -2\, \tilde{R}^\eps_x \wedge \tilde{R}^\eps_y +h^\eps ,
 \ee
where $h^\eps\rightarrow 0$ in $L^1_{loc}(\R^2) $ as $\eps \rightarrow 0$. Then, up to a subsequence, we have
\be
\tilde{R}^\eps \rightarrow R  \hbox{ a.e. on } \R^2
\ee
and
\be
\nabla R^\eps \rightharpoonup \nabla R \hbox{ weakly in } L^2(\R^2).
\ee
Moreover $R$ is a weak solution of (\ref{eqlim}). Thanks to our choice of $\nu$, we are going to prove that the weak convergence is in fact a strong convergence. Let $v^\eps = \tilde{R}^\eps -R$, then $v^\eps$ satisfies
\be
\Delta v^\eps= -2\, v^\eps_x \wedge v^\eps_y -2 (v^\eps_x \wedge R_y + R_x \wedge v^\eps_y) +h^\eps .
\ee
Thanks to lemma \ref{wente2}, there exists $\psi_\eps$ a the solution in $H^1(\R^2)$ of 
\be
\Delta \psi^\eps = -2 (v^\eps_x\wedge R_y+ R_x \wedge v^\eps_y),
\ee
which satisfies
\beq
\label{h1}
\Vert \nabla \psi^\eps  \Vert_2 + \Vert \psi^\eps  \Vert_\infty  \leq \Vert \nabla v^\eps  \Vert_2 \Vert \nabla R  \Vert_2  .
\eeq
On the other hand,
\be
\int_{\R^2} \vert \nabla \psi^\eps \vert^2 dz = -2 \int_{\R^2}  \langle \psi^\eps , v^\eps_x\wedge R_y+ R_x \wedge v^\eps_y\rangle dz .
\ee
Then, thanks to (\ref{h1}),  $\psi^\eps \wedge R_x$ and  $\psi^\eps \wedge R_y$  are bounded in $L^2(\R^2)$. Hence, since $\nabla v^\eps\rightarrow 0$ weakly in $L^2$, it follows that 
$$\int_{\R^2} \vert \nabla \psi^\eps \vert^2 dz \rightarrow 0. $$
Finally we have
$$ \Delta v^\eps = -2\, v^\eps_x \wedge v^\eps_y  + g^\eps , $$
where $g^\eps\rightarrow 0$ in $D'(\R^2)$.\\ 

Finally, let $\phi \in C^\infty_c(\R^2)$ such that $supp(\phi)$ is contained in a ball of radius $1$, using lemma \ref{wente3}, we have
\be
\begin{split}
\int_{\R^2} \vert \nabla (\phi v^\eps) \vert^2 dz &= - 2 \int_{\R^2}  \langle  v^\eps , \phi  v^\eps_x\wedge  \phi v^\eps_y \rangle dz  + o(1) ,\\
& \leq 2\left(C_0\Vert \nabla v^\eps_{| supp(\phi)} \Vert_2 \right) \Vert \nabla(\phi  v^\eps)\Vert^2_2 +o(1). 
\end{split}
\ee
Thanks to our choice of $\lambda^\eps$, we have $C_0 \Vert \nabla v^\eps_{| supp(\phi)} \Vert_2 \leq \frac{1}{2}$, and then 
$$ \int_{\R^2} \vert \nabla (\phi v^\eps) \vert^2 dz = o(1)$$
which prove that 
$$\nabla \tilde{R}^\eps \rightarrow \nabla R \hbox{ strongly in } L^2_{loc}(\R^2) .$$
Then the convergence is strong, since it was already the case far from $0$. Moreover $R$ is not constant since $\Vert \nabla R \Vert_2 =\nu>0$. But thanks to (\ref{ii}) we get that, for all $z\in \R^2$, there exists $i$ such that 
$$ \vert \nabla \tilde{R}^\eps (z)\vert =o\left( \frac{1}{\sqrt{ \left(\frac{\lambda_{i}^\eps}{\lambda^\eps}\right)^2 + \left\vert z +\frac{a^\eps -a^\eps_i}{\lambda^\eps}\right\vert^2} }\right) ,$$  
 which leads to a contradiction and proves (\ref{h1dec}).\\
 
Finally, in order to finish the proof of the theorem we just need to prove that $max\{ deg P_i, deg\ Q_i\}= 1$ for all $0 \leq i\leq N$, with  $\omega_i =\pi_{P_i}\left( \frac{P_i}{Q_i}\right) $ where $\frac{P_i}{Q_i}$ is irreducible. But this is an easy consequence of the fact that our surfaces are embedded, (\ref{A}) and  lemma \ref{simple} and this achieves the proof of the theorem.\hfill$\square$ \\

%%%%%%%%%%%%%%%%%%%%%%%%%%%%%%%%%%%%%%%%%%%%%%%%%%%%%%
In the previous theorem, we showed that $u^\eps$ behaves asymptotically as a sum of Euclidean bubbles. Hence in this decomposition the curved term doesn't play any role. In order to prove our theorem, we have to make appear these terms and to be more precise about this decomposition. In this goal,  we expand the metric in (\ref{main3}) thanks to appendix \ref{a4}, which gives
\beq
\label{exp}
\begin{split}
\Delta u^\eps_j &= - 2\; (u_x^\eps \wedge u_y^\eps)_j + \eps^2 \left( \frac{2}{3} R_{imnj}(p_\eps) (u^\eps)^m (u^\eps)^n (u_x^\eps \wedge  u_y^\eps )^i\right. \\
&+\frac{1}{3}\mathrm{Ric}_{mn}(p_\eps)(u^\eps)^m (u^\eps)^n(u_x^\eps \wedge  u_y^\eps )_j \\
&\left. +\frac{1}{3}\left( R_{nmij}(p_\eps) + R_{njim}(p_\eps)\right)(u^\eps)^m \langle \nabla (u^\eps)^i,  \nabla (u^\eps)^n \rangle  \right) + O(\eps^3 \vert \nabla u_\eps \vert^2),
\end{split}
\eeq
where $p_\eps$ is the center of our chart.\\

Since we want to prove something about the derivatives of the curvature, we have to eliminate the curvature term. In order to do it we are going to be more precise on  the shape of our bubbles. In fact this bubbles aren't euclidian once $\eps>0$, the curvature make them look like ellipsoids. For each $i$ we are going to search a perturbation of $\omega^i$ such that when we consider (\ref{exp}) around $a_i^\eps$ the curvature term disappear. A last thing we have to pay attention is that our $\omega^i$ are {\it a priori} not centered at $0$, which is not very convenient when we want to compare $\omega^i_x \wedge \omega^i_y$ and $\omega^i$. Then we will consider sometimes $\hat{\omega}^i =\omega^i -p^i$ where $p^i$ is the center of mass of $\omega^i$.\\
 
 Hence we look for $\rho_i^\eps$ such that $\overline{\omega}_i^\eps= \hat{\omega}^i + p_i^\eps+ \rho^\eps_i$ solves (\ref{exp}) at the first order, here $p_i^\eps$ is a constant we will fix later. That is to say
\beq
\begin{split}
\Delta \overline{\omega}_i^\eps &=  -2 (\overline{\omega}_i^\eps)_x \wedge (\overline{\omega}_i^\eps)_y\\
 &+ \eps^2 \left( \frac{2}{3} R_{kmnl}(p_\eps) (\overline{\omega}_i^\eps)^m (\overline{\omega}_i^\eps)^n ((\overline{\omega}_i^\eps)_x \wedge  (\overline{\omega}_i^\eps)_y)^k\right. \\
 &+ \frac{1}{3}\mathrm{Ric}_{mn}(p_\eps)(\overline{\omega}_i^\eps)^m (\overline{\omega}_i^\eps)^n((\overline{\omega}_i^\eps)_x \wedge  (\overline{\omega}_i^\eps)_y)_l\\ 
  &\left. + \frac{1}{3}\left( R_{mnkl}(p_\eps) + R_{mlkn}(p_\eps)\right)(\overline{\omega}_i^\eps)^n \langle  \nabla (\overline{\omega}_i^\eps)^k,  \nabla (\overline{\omega}_i^\eps)^m \rangle  \right) + O(\eps^3)
\end{split}
\eeq
with the relation of almost conformality
\beq
\label{conforme}
\begin{split}
& \langle( \overline{\omega}_i^\eps)_x,  ( \overline{\omega}^\eps_i)_y\rangle +\frac{1}{3} R_{kmnl}(p_\eps) (\overline{\omega}_i^\eps)^m (\overline{\omega}_i^\eps)^n  (\overline{\omega}_i^\eps)_x^{k}(\overline{\omega}_i^\eps)_y^{l} =O(\eps^3) ,\\
& \langle(\overline{\omega}_i^\eps)_x, ( \overline{\omega}_i^\eps)_x\rangle +\frac{1}{6} R_{kmnl}(p_\eps) (\overline{\omega}_i^\eps)^m (\overline{\omega}_i^\eps)^n  (\overline{\omega}_i^\eps)_x^{k}(\overline{\omega}_i^\eps)_x^{l} \\
&- \langle (\overline{\omega}_i^\eps)_y, (\overline{\omega}_i^\eps)_y\rangle - \frac{1}{6} R_{kmnl}(p_\eps) (\overline{\omega}_i^\eps)^m (\overline{\omega}_i^\eps)^n  (\overline{\omega}_i^\eps)_y^{k}(\overline{\omega}_i^\eps)_y^{l} =O(\eps^3). 
\end{split}
\eeq
 We look for $\rho^\eps_i$ of the form  $\eps^2 \rho^\eps_i$. Then, thanks to the expansion of the metric and~(\ref{di}), then we see that  $\rho_i^\eps$ must solve
\beq
\begin{split}
&\Delta \rho^\eps_i +2 ((\rho^\eps_i)_x\wedge \hat{\omega}^i_y + \hat{\omega}^i_x \wedge (\rho^\eps_i)_y ) =  \left( -\frac{2}{3} R_{kmnl}(p_\eps) (\hat{\omega}_i +p_i^\eps )^m (\hat{\omega}_i +p_i^\eps)^n (\hat{\omega}_i)^k \right. \\
 &- \frac{1}{3}\mathrm{Ric}_{mn}(p_\eps)(\hat{\omega}_i +p_i^\eps)^m (\hat{\omega}_i +p_i^\eps)^n(\hat{\omega}_i)_l\\ 
  &\left. + \left(\frac{1}{3} R_{mnkl}(p_\eps) + R_{mlkn}(p_\eps)\right)(\hat{\omega}_i +p_i^\eps)^n (\delta_{km} -  \hat{\omega}_i^k   \hat{\omega}_i^m)  \right) \frac{\vert \nabla \hat{\omega}_i \vert^2}{2}
 \end{split}
 \eeq
and
\be
\begin{split}
\langle (\rho^\eps_i)_x, \hat{\omega}^i_y \rangle + \langle \hat{\omega}^i_x, (\rho^\eps_i)_y \rangle &= -\frac{1}{3} R_{kmnl}(p_\eps)(\hat{\omega}_i +p_i^\eps)^m (\hat{\omega}_i +p_i^\eps)^n (\hat{\omega}^i_x)^{k}(\hat{\omega}^i_y)^{l},\\
\langle (\rho^\eps_i)_x, \hat{\omega}^i_x \rangle - \langle \hat{\omega}^i_y, (\rho^\eps_i)_y \rangle &=  \frac{1}{6} R_{kmnl}(p_\eps)(\hat{\omega}_i +p_i^\eps)^m (\hat{\omega}_i +p_i^\eps)^n (\hat{\omega}^i_y)^{k}(\hat{\omega}^i_y)^{l}\\
& -\frac{1}{6} R_{kmnl}(p_\eps) (\hat{\omega}_i +p_i^\eps)^m (\hat{\omega}_i +p_i^\eps)^n (\hat{\omega}^i_x)^{k}(\hat{\omega}^i_x)^{l}.
\end{split}
\ee
As for the linearized equation, see  proposition \ref{lin}, we decompose $\nabla \rho^\eps_i$ on the orthogonal frame $\hat{\omega}^i_x, \hat{\omega}^i_y, \hat{\omega}^i_x\wedge \hat{\omega}^i_y$ in order to find the solution. After a straitforward computation, we check that  
\beq
\label{formula}
\begin{split}
\rho_i^\eps &= \frac{1}{6} \left(\mathrm{Ric}_{kl}(p_\eps) (\hat{\omega}^i)^l  - \frac{3}{2}Scal(p^\eps) (\hat{\omega}^i)^k\right)- \frac{1}{6} R_{kmnl}(p_\eps) (p_i^\eps)^m (p_i^\eps)^n (\hat{\omega}^i)^l \\ 
&- \frac{1}{3} R_{kmnl}(p_\eps) (p_i^\eps)^m (\hat{\omega}^i)^n (\hat{\omega}^i)^k - \frac{1}{12} \hbox{Ric}_{kl}(p_\eps) (\hat{\omega}^i)^k (\hat{\omega}^i)^l  \hat{\omega}^i 
\end{split}
\eeq  
provides a solution. Here we used the fact that in dimension 3, we get 
\be
\begin{split}
R_{kmnl} &=(g_{kn}\mathrm{Ric}_{ml}-g_{kl}\mathrm{Ric}_{mn} + g_{ml}\mathrm{Ric}_{kn} - g_{mn}\mathrm{Ric}_{kl}) \\
&+\frac{\hbox{Scal}}{2} (g_{kl}g_{mn}-g_{kn}g_{ml}).
\end{split}
\ee
Hence we  set 
\beq
\label{defb}
\begin{split}
&\overline{\omega}_i^\eps= \hat{\omega}^i + p_i^\eps + \eps^2 \rho_i^\eps \\
&\hbox{and }\\
&B_i^\eps(z) =\overline{\omega}_i^\eps \left( \frac{z-a^\eps_i}{\lambda_i^\eps} \right).
\end{split}
\eeq
But, we have to make an adjustment on our bubbles, choosing them such that they are  "tangent" to $\Sigma_\eps$ at its extreme points.\\

 First, let $i\in \{0,\dots, p\}$ be such that 
\beq
\label{hyp0}
\lim_{\eps \rightarrow 0} \frac{d_j^\eps(a_i^\eps)}{\lambda_i^\eps} =0 \hbox{ for some  } j\not=i .
\eeq
Let us fix $b_i^\eps \in R^2$ such that 
\be
\begin{cases}
b_i^\eps \in B(a_i^\eps, \lambda_i^\eps)\\
\hbox{ and }\\
d\left(\frac{b_i^\eps -a_i^\eps}{\lambda_i^\eps}, S_i \right) \geq d>0 ,
 \end{cases}
 \ee
and $p_i^\eps\in \R^3$ such that 
\be
p_i^\eps + \hat{\omega}^i(b_i^\eps)=u^\eps(b_i^\eps). 
 \ee
Then, we consider $i\in \{0,\dots, p\}$  such that 
\beq
\label{hyp}
\lim_{\eps \rightarrow 0} \frac{d_j^\eps(a_i^\eps)}{\lambda_i^\eps} \not=0 \hbox{ for any } j\not=i .
\eeq
Hence, thanks to theorem \ref{de}, there exists $\delta_0>0$ such that, up to a subsequence,
\beq
\label{conver}
\nabla \tilde{u}^\eps  \rightarrow \nabla \omega^{i} \hbox { in } C^2(B(0,\delta_0)) ,
 \eeq
where  $\tilde{u}^\eps = u^\eps(\lambda_{i}^\eps\, . + a_{i}^\eps ) $. In fact the convergence should hold on $B(0, \delta_0)\setminus \{S_{i}\}$. But, thanks to (\ref{hyp}) and (\ref{ortho}), either  $\ds a_{ij} =\lim_{\eps\rightarrow 0}\frac{a_{j}^\eps - a_{i}^\eps  }{\lambda_{i}^\eps}\not =0$ or $\lambda_{i}^\eps =o(\lambda_j^\eps)$. Hence, in every cases, we get that 
$$ \vert \nabla\tilde{B}_j^\eps \vert \rightarrow 0 \hbox { in } C^1(B(0,\delta_0)), \hbox{ for all } j\not= i $$
which proves the validity of (\ref{conver}).\\

Then, thanks to (\ref{conver}) and the fact that $\vert\nabla \omega^{i}\vert$ has a strict maximum at $0$, for $\eps$ small enough, there exists $\tilde{a}^\eps_{i} \in\R^2$ such 
\beq
\label{xx}
\vert \tilde{a}^\eps_{i}- a^\eps_{i}\vert =o(\lambda_{i}^\eps) \hbox{ and } \left\vert \nabla u^\eps \right\vert \hbox{ has a local maximum at } \tilde{a}^\eps_{i} .
\eeq
Still thanks to (\ref{conver}), there exists $R_{i}^\eps \in SO(3)$, $\theta^\eps_i \in [0,2\pi]$ and $\tilde{\lambda}_i^\eps\in \R$ such that
\beq
\label{yy}
\begin{split}
&\tilde{\lambda}_i^\eps \sim\lambda_i^\eps,\; R_i^\eps \rightarrow Id ,\\
&\tilde{u}^\eps_x (0)= R_i^\eps (\omega^i_x (0))\\
&\hbox{ and } \\
&\hbox{Vect}\left( \tilde{u}^\eps_x (0), \tilde{u}^\eps_y (0)\right)= \hbox{Vect}\left( R_i^\eps (\omega^i_x (0)),R_i^\eps (\omega^i_y (0))\right),
\end{split}
\eeq
where $\ds \tilde{u}^\eps=  u^\eps (e^{i\theta_i^\eps}\tilde{\lambda}_{i}^\eps\, .\,+\tilde{a}^\eps_{i})$. Then, we set $p_i^\eps\in \R^3$ such that 
\be
\begin{split}
\label{pp}
&\hat{\omega}_i^\eps = R_i^\eps \hat{\omega}^i ,\\
&p_i^\eps=u^\eps(\tilde{a}_i^\eps) - \hat{\omega}_i^\eps (0) \\
&\overline{\omega}_i^\eps= \hat{\omega}_i^\eps + p_i^\eps + \eps^2 \rho_i^\eps , \\
&\hbox{and }\\
&B_i^\eps(z) = \overline{\omega}_i^\eps \left( \frac{z-\tilde{a}_i^\eps}{e^{i\theta_i^\eps}\tilde{\lambda}_i^\eps} \right) ,
\end{split}
\ee
where $\rho_i^\eps $ is associated to $\hat{\omega}^i_\eps$  and $p_i^\eps$ thanks to (\ref{formula}).\\

Moreover, thanks to (\ref{conver}), there exists $c_i^\eps\in\R^2$ such that
\beq
\label{ww}
\vert c^\eps_{i} \vert =o(1) \hbox{ and } \left\vert \nabla \overline{\omega}_i^\eps \right\vert \hbox{ has a local maximum at } c^\eps_{i} .
\eeq
However using the fact that $\vert \nabla \hat{\omega}^\eps_i\vert $ has a local maximum at $0$ and the fact that  $\vert \nabla (\hat{\omega}^\eps_i-\overline{\omega}_i^\eps)\vert=O(\eps^2) $ in a neighborhood of $0$, we get that 
\be
\begin{split}
& \vert c^\eps_{i} \vert =O(\eps)\\
&\hbox{ and }\\
& \vert \nabla \tilde{u}^\eps (0) -\nabla  \overline{\omega}_i^\eps (c^\eps_i)\vert =O(\eps)
\end{split}
\ee
where $\ds \tilde{u}^\eps=  u^\eps (e^{i\theta_i^\eps}\tilde{\lambda}_{i}^\eps\, .\,+\tilde{a}^\eps_{i})$. This implies that there exist $\hat{R}_{i}^\eps \in SO(3)$  $\hat{\theta}^\eps_i \in [0,2\pi]$ and $\hat{\lambda}_{i}^\eps\in \R$ such that
\beq
\label{zz}
\begin{split}
&\left\vert \frac{\hat{\lambda}_i^\eps}{\lambda_i^\eps}- 1\right\vert=O(\eps),\; \vert \hat{R}_i^\eps - Id\vert =O(\eps) ,\\
&u^\eps_x (\tilde{a}_i^\eps)= \hat{R}_i^\eps ((\hat{B}^\eps_i)_x (\tilde{a}_i^\eps))\\
&\hbox{ and } \\
&\hbox{Vect} \left( u^\eps_x (\tilde{a}_i^\eps), u^\eps_y (\tilde{a}_i^\eps)\right)= \hbox{Vect}\left( \hat{R}_i^\eps ((\hat{B}^\eps_i)_x (\tilde{a}_i^\eps)), \hat{R}_i^\eps ((\hat{B}^\eps_i)_y (\tilde{a}_i^\eps))\right),
\end{split}
\eeq
where $\hat{B}^\eps_i= \overline{\omega}_i^\eps \left( \frac{z-\tilde{a}_i^\eps}{e^{i\hat{\theta}_i^\eps}\hat{\lambda}_i^\eps} + c^\eps_i \right)$.
Then we replace $a^i_\eps$ by $\tilde{a}_i^\eps$, $\lambda_i^\eps$ by $ e^{i\theta_i^\eps}\hat{\lambda}_{i}^\eps$ and   $B_i^\eps$ by $\hat{R}_i^\eps \hat{B}_i^\eps$. Thanks to (\ref{xx}), (\ref{yy}) and  (\ref{zz}) the conclusions of theorem \ref{de} still holds with our new choice of $\lambda_i^\eps $ and $a_i^\eps$. Moreover, thanks to(\ref{conforme}), (\ref{pp}), (\ref{zz}) and the fact that we have adjust the tangent plane, for every $i$ that satisfies (\ref{hyp}), we have 
\beq
\label{ladjust}
\begin{split}
&\left(  \tilde{u}^\eps -\tilde{B}_{i}^\eps \right)(0) = O\left(\eps \right),\\
&\left\vert \nabla\left( \tilde{u}_i^\eps - \tilde{B}_{i}^\eps \right)(0)\right\vert^2 = O\left(\eps^3\right),\\
& \hbox{ and }\\
&\left\vert \nabla^2\left( \tilde{u}_i^\eps - \tilde{B}_{i}^\eps \right)\left( \nabla \tilde{B}_i^\eps\right)(0)\right\vert^2 = O\left(\eps^3\right),
\end{split}
\eeq
where $\ds \tilde{f}^\eps_i= f(\lambda_{i}^\eps\, .\,+ a^\eps_{i})$. Then, we give the equations satisfied by our modified bubbles,
\beq
\label{eqb}
\begin{split}
\Delta B_i^\eps &= -2 (B_i^\eps)_x \wedge (B_i^\eps)_y +\eps^2 \left( \frac{2}{3} R_{imnj}(p_\eps) (B_i^\eps)^m (B_i^\eps)^n ((B_i^\eps)_x \wedge  (B_i^\eps)_y )^i\right. \\
&+\frac{1}{3}\mathrm{Ric}_{mn}(p_\eps)(B_i^\eps)^m (B_i^\eps)^n ((B_i^\eps)_x \wedge  (B_i^\eps)_y )_j \\
&\left. +\frac{1}{3}\left( R_{mnij}(p_\eps) + R_{mjin}(p_\eps)\right)(B_i^\eps)^n \langle \nabla (B_i^\eps)^i,  \nabla (B_i^\eps)^m \rangle  \right) +  O(\eps^3 \vert \nabla B_i^\eps\vert^2)
\end{split}
\eeq
and the relation of quasi-conformality
\beq
\label{eqbc}
\begin{split}
&\langle (B_i^\eps)_x, (B_i^\eps)_y \rangle = -\frac{1}{3} R_{kmnl}(p_\eps) (B_i^\eps)^m (B_i^\eps)^n  (B_i^\eps)_x^{k}(B_i^\eps)_y^{l}+ O(\eps^3 \vert \nabla B_i^\eps\vert^2),\\
&\langle (B_i^\eps)_y, (B_i^\eps)_y \rangle- \langle (B_i^\eps)_x, (B_i^\eps)_x \rangle = \frac{1}{6} R_{kmnl}(p_\eps) (B_i^\eps)^m (B_i^\eps)^n  (B_i^\eps)_y^{k}(B_i^\eps)_y^{l} \\
&-\frac{1}{6} R_{kmnl}(p_\eps) (B_i^\eps)^m (B_i^\eps)^n  (B_i^\eps)_x^{k}(B_i^\eps)_x^{l}+ O(\eps^3 \vert \nabla B_i^\eps\vert^2).
\end{split}
\eeq
{\bf Conclusion of the decomposition step :}\\

Finally, let $u^\eps$ be a sequence of $C^2$-solutions of (\ref{main3}). Then, there exist $p \in \N$ and 
\begin{enumerate}[(i)]
\item   $\omega^0,\dots,\omega^p $ simple solutions of (\ref{eqlim}) such that $\vert \nabla \omega^i \vert$ has a maximum at $0$, 
\item  $ a^\eps_0, \dots, a^\eps_p$ bounded sequences of $\R^2$ such that  $\ds \lim_{\eps\rightarrow 0} a_i^\eps = 0$ for all $1\leq i \leq p$, 
\item  $\lambda^\eps_0 , \dots , \lambda^\eps_p$ bounded sequences of complex numbers such that $\ds \lim_{\eps \rightarrow 0} \lambda^\eps_i =0$ for all $1\leq i \leq p$, 
\item $ R^\eps_0, \dots, R^\eps_{2p+1}$  sequences of $SO(3)$ such that  $\ds\lim_{\eps\rightarrow 0} R_i^\eps = Id$ for all $0\leq i \leq 2p+1$,
\item $ p^\eps_0, \dots, p^\eps_p$  sequences of points of $\R^3$ such that  $\ds\lim_{\eps\rightarrow 0} p_i^\eps = p_i$ for all $0\leq i \leq p$, where $p_i$ is the center of mass of $\omega^i$, 
\end{enumerate} 
such that, for a subsequence of $u_\eps$ (still denoted $u^\eps$) the following assertions hold 
\be
 u^\eps_i \rightarrow \omega^i \hbox{ in } C^2_{loc}(\R^2\setminus S_i) \hbox{ as } \eps \rightarrow 0 \hbox{ for all } 0\leq i\leq p,
 \ee
 where $u^\eps_i= u^\eps (\lambda^\eps_i\ . + a_i^\eps )$ and $\ds S_i= \lim_{\eps \rightarrow 0} \left\{ \frac{a_j^\eps-a_i^\eps}{\lambda_i^\eps}\hbox{ s.t. } j\in \{0,\dots,p\}\setminus \{ i\}\right\} $.
\beq
\label{orthof}
 \frac{d_{i}^\eps(a_j^\eps)}{\lambda^\eps_j} + \frac{d_{j}^\eps(a_i^\eps)}{\lambda^\eps_i} \rightarrow +\infty  \hbox{ for all }i \not= j ,
\eeq
 where $d_{i}^\eps(x)=\sqrt{(\lambda^\eps_i)^2 +\vert a_i^\eps -x\vert^2}$.\\
 
With the additional properties that 
\beq
\label{estfaible}
 \lim_{\eps \rightarrow 0} \sup_{x\in \R^2} \left(\min_{0 \leq i \leq p} d_i^\eps (x)\right) \left\vert \nabla  \left(u^\eps-\sum_{i=0}^p  B^\eps_i\right)  (x) \right\vert =0
\eeq
and 
\be
\left\Vert\nabla \left( u^\eps -\sum_{i=0}^p B^\eps_i \right) \right\Vert_2 \rightarrow 0  \hbox{ as } \eps\rightarrow 0,
\ee 
where 
\be
\begin{split}
& \hat{\omega}^i  = \omega^i -p_i,\\
&\hat{\omega}_i^\eps =  R_i^\eps \hat{\omega}^i,\\
&\overline{\omega}_i^\eps= \hat{\omega}_i^\eps + p_i^\eps + \eps^2 \rho_i^\eps , \\
&\hbox{and }\\
&B_i^\eps(z) = R_{p+1+i}^\eps \overline{\omega}_i^\eps \left( \frac{z-a_i^\eps}{\lambda_i^\eps} \right) ,
\end{split}
\ee
where $\rho_i^\eps $ is associated to $\hat{\omega}^i_\eps$ and $p^\eps_i$ thanks to (\ref{formula}).\\

Moreover, for all $i$,  there exists $b_i^\eps \in  R^2$ such that 
\beq
\label{ffa}
\begin{split}
&b_i^\eps \in B(a_i^\eps, \lambda_i^\eps)\\
&d\left(\frac{b_i^\eps -a_i^\eps}{\lambda_i^\eps}, S_i \right) \geq d>0 ,\\
&\hbox{ and }\\
&\vert u^\eps(b^\eps_i)- B_i^\eps(b^\eps_i)\vert =O(\eps),
 \end{split}
 \eeq
and for all $i$ such that  $\ds  \lim_{\eps \rightarrow 0} \frac{d_j^\eps(a_i^\eps)}{\lambda_i^\eps} \not=0 \hbox{ for any } j\not=i$, we get
\beq
\label{cint}
\begin{split}
&\left\vert \nabla\left( \tilde{u}^\eps -\tilde{B}_{i}^\eps \right)(0)\right\vert^2 = O\left(\eps^3\right),\\
&\left\vert \nabla^2 \left(  \tilde{u}^\eps -\tilde{B}_{i}^\eps \right)(\tilde{B}_{i}^\eps)(0)\right\vert^2 = O\left(\eps^3\right).
\end{split}
\eeq
Finally, we can replace the $\lambda_i^\eps$ by there absolute value since the argument part can be in the $B^\eps_i$. Moreover, thanks to the method we used to construct our $\lambda_i^\eps$, up to reordering, we can assume that 
$$\Vert \nabla u^\eps\Vert_{\infty} =\frac{1}{\lambda_p^\eps}.$$  

%%%%%%%%%%%%%%%%%%%%%%%%%%%%%%
\section{Strong estimate}
Let $B_i^\eps$ defined as in the last section and  $\ds R^\eps = u^\eps -\sum_{i=0}^p B_i^\eps$ be the remainder. The aim of this step is to prove an estimate on the gradient of the remainder, $r^\eps=\Vert \nabla R^\eps \Vert_{\infty}$. Thanks to the previous step  $R^\eps$ satisfies the following equations
\beq
\label{eqr}
\begin{split}
\Delta R^\eps &= -2 \left(  \sum_{i\not= j} (B_i^\eps)_x \wedge (B_j^\eps)_y + \sum_{i=0}^p  (B_i^\eps)_x \wedge R^\eps_y + R^\eps_x \wedge (B_i^\eps)_y \right)\\
&-2\, R^\eps_x \wedge R^\eps_y +O\left(\eps^2 \left( \sum_{i=0}^p  \vert u^\eps - B^\eps_i \vert \vert \nabla B_i^\eps\vert^2 \right) \right)\\
&+ O\left(\eps^2 \left(  \sum_{i=0}^p  \vert \nabla B_i^\eps\vert \left( \sum_{j\not= i} \vert \nabla B_j^\eps\vert +\vert \nabla R^\eps\vert\right)  + \vert \nabla R^\eps\vert^2  \right)\right)\\
&+\left( \eps^3 \vert \nabla u^\eps\vert^2 \right) ,
\end{split}
\eeq
\beq
\label{eqrc}
\begin{split}
&\sum_{i\not= j} \langle (B_i^\eps)_x , (B_j^\eps)_y \rangle + \langle R^\eps_x ,R^\eps_y\rangle + \sum_{i=0}^p  \langle (B_i^\eps)_x , R^\eps_y \rangle + \langle R^\eps_x , (B_i^\eps)_y \rangle \\ 
&=  O\left(\eps^2 \left( \sum_{i=0}^p  \vert u^\eps - B^\eps_i \vert \vert \nabla B_i^\eps\vert^2 \right) \right)+\left( \eps^3 \vert \nabla u^\eps\vert^2 \right)\\
&+ O\left(\eps^2 \left(  \sum_{i=0}^p  \vert \nabla B_i^\eps\vert \left( \sum_{j\not= i} \vert \nabla B_j^\eps\vert +\vert \nabla R^\eps\vert\right)  + \vert \nabla R^\eps\vert^2  \right)\right),
\hbox{ and }\\
&\sum_{i\not= j} \langle (B_i^\eps)_x , (B_j^\eps)_x \rangle- \langle (B_i^\eps)_y , (B_j^\eps)_y \rangle+2 \sum_{i=0}^p  \langle (B_i^\eps)_x , R^\eps_x \rangle -\langle  (B_i^\eps)_y, R^\eps_y \rangle \\
& + \langle R^\eps_x ,R^\eps_x\rangle - \langle R^\eps_y ,R^\eps_y\rangle =  O\left(\eps^2 \left( \sum_{i=0}^p  \vert u^\eps - B^\eps_i \vert \vert \nabla B_i^\eps\vert^2 \right) \right)+\left( \eps^3 \vert \nabla u^\eps\vert^2 \right)\\
&+ O\left(\eps^2 \left(  \sum_{i=0}^p  \vert \nabla B_i^\eps\vert \left( \sum_{j\not= i} \vert \nabla B_j^\eps\vert +\vert \nabla R^\eps\vert\right)  + \vert \nabla R^\eps\vert^2  \right)\right)\end{split}
\eeq
Then our aim  is to show that the remainder is controlled  by $\eps^3 \Vert \nabla u^\eps\Vert_\infty=\frac{\eps^3}{\lambda_p^\eps}$. Indeed if the contrary holds then the last terms of (\ref{eqr}) and (\ref{eqrc}) would be negligible and we would get a non-trivial solutions of the linearized problem which has only trivial solution for a good choice of the initial data, see proposition \ref{lin}.\\

But  we also need to control the cross terms like  $\ds \sum_{j\not= i}  \vert \nabla B_i^\eps\vert \vert \nabla B_j^\eps\vert$. Hence, we define, for $i\not=j$, 
$$t_{ij}^\eps = \frac{\lambda_j^\eps}{(d_j^\eps(a_i^\eps))^2 +(d_i^\eps(a_j^\eps))^2}$$
which is the trace around $a_i^\eps$ of $\nabla B^\eps_j$. Indeed, thanks to (\ref{oi}), we easily check that on any compact subset of  $ \R^2\setminus{\{ S_i\}}$, we have
$$\left\vert \nabla B_j^\eps \left(\lambda_i^\eps\, . + a_i^\eps \right) \right\vert \sim c\, t^\eps_{ij},$$
where $c$ is a positive constant. Then we define also the maximum of these interactions as
$$t^\eps =\max_{i\not= j}(t^\eps_{ij}).$$
In order to get an estimate on $R^\eps$, the idea is to apply Green identity to (\ref{eqr}). This is possible thanks to lemma \ref{green}. Hence let $z^\eps\in\R^2$, we get 
\beq
 \label{green2}
 \begin{split}
\vert \nabla R^\eps (z^\eps) \vert  &\leq I_{i}^\eps(z^\eps) +I_{ij}^\eps(z^\eps) + \vert \nabla \phi^\eps(z^\eps)\vert \\
 &+O\left( \eps^2 \left( \sum_{i=0}^{p} (J_{i}^\eps(z^\eps) +I_{i}^\eps(z^\eps)) + \sum_{0\leq i<j\leq p } I_{ij}^\eps(z^\eps)   \right)  \right) +O(I^\eps(z^\eps))\\
 & +O( J^\eps(z^\eps))
 \end{split}
\eeq
where 
\be
\begin{split}
I_i^\eps(z^\eps) &= \int_{\R^2} \vert \nabla G(\, .\, ,z^\eps) \vert \vert \nabla B_i^\eps \vert\vert \nabla R^\eps \vert dz \\
I_{ij}^\eps(z^\eps) &= \int_{\R^2} \vert \nabla G(\, .\,,z^\eps) \vert \vert \nabla B_i^\eps \vert \vert \nabla B_j^\eps \vert dz\\
J_i^\eps(z^\eps) &=   \int_{\R^2}  \vert \nabla G(\, .\,,z^\eps) \vert \vert u^\eps - B^\eps_i \vert \vert \nabla B_i^\eps\vert^2 dz \\
I^\eps(z^\eps) &=   \eps^2 \int_{\R^2} \vert \nabla G(\, .\,,z^\eps) \vert \vert \nabla R^\eps \vert^2 dz \\
J^\eps(z^\eps) &=  \eps^3 \int_{\R^2} \vert \nabla G(\, .\,,z^\eps) \vert \vert \nabla u^\eps \vert^2 dz 
 \end{split}
\ee
and $\phi^\eps\in H^1(\R^2)$ est une solution de
$$\Delta\phi^\eps = -2  R^\eps_x \wedge R^\eps_y .$$  
First, using lemma \ref{wente2} and  (\ref{h1dec}), we deduce that
\beq
\label{xxr2} 
\vert \nabla \phi^\eps(z^\eps)\vert =O (\Vert \nabla R^\eps\Vert_2 \Vert \nabla R^\eps\Vert_\infty)=o(r^\eps) .
\eeq
Then, we estimate $ I_{i}^\eps(z^\eps)$, $I_{ij}^\eps(z^\eps)$ and $J_i^\eps(z^\eps)$.\\
Let $R>0$, we set $\ds r^\eps _{i,R}= \sup_{\Omega_{i,R}^\eps}\vert \nabla R^\eps \vert$ where $\ds\Omega_{i,R}^\eps = B(a_{i}^\eps, \lambda_{i}^\eps R)\setminus \left\{\cup_{j\not= i} B(a_j^\eps, \frac{\lambda_{i}^\eps}{R} ) \right\} $. For $\eps$ small enough, we get
\beq
\label{Ij2}
\begin{split}
 I_{i}^\eps(z^\eps) &\leq r^\eps \int_{\R^2 \setminus \Omega_{i,R}^\eps } \vert \nabla G(\, .\,,z^\eps) \vert \vert \nabla B_i^\eps  \vert dz\\
  &+ r_{i,R}^\eps \int_{\Omega_{i,R}^\eps } \vert \nabla G(\, .\,,z^\eps) \vert \vert \nabla B_i^\eps  \vert dz .
  \end{split}
  \eeq
  But, by a simple change of variable we get,  
  $$ \int_\Omega \vert \nabla G(\, .\,,z^\eps) \vert \vert \nabla B_i^\eps  \vert dz  = O\left( \int_{\frac{\Omega -a_i^\eps}{\lambda_i^\eps}}  \left\vert \nabla G \left(\, .\, ,\frac{z^\eps -a_i^\eps}{\lambda_i^\eps} \right) \right\vert \frac{1}{1+\vert z\vert^2} \vert dz\right),$$
where $\Omega$ is any measurable set. Then either $ \frac{z^\eps -a_i^\eps}{\lambda_i^\eps} \rightarrow +\infty$ and ,thanks to lemma \ref{estimint1}, we get 
  $$ \int_{\R^2} \vert \nabla G(\, .\,,z^\eps) \vert \vert \nabla B_i^\eps  \vert dz =o(1)$$
or  $ \frac{z^\eps -a_i^\eps}{\lambda_i^\eps} \rightarrow z_0$ and 
$$ \lim_{\eps\rightarrow 0}\int_{\R^2 \setminus \Omega_{i,R}^\eps } \vert \nabla G(\, .\,,z^\eps) \vert \vert \nabla B_i^\eps  \vert dz =O\left( \int_{\R^2 \setminus \Omega_{i,R}}  \vert \nabla G(\, .\,,z_0) \vert \frac{1}{1+\vert z\vert^2} \vert dz\right) $$
where $\Omega_{i,R} = B(0, R)\setminus \left\{\cup_{z\in S_{i}} B(z, \frac{1}{R} ) \right\} $. Hence in every case we get 
$$
 I^\eps_i (z) \leq r^\eps \delta_{R,\eps} +  O\left(r_{i,R}^\eps \frac{Ln\left( 2+ \frac{\vert a_i^\eps -z^\eps\vert}{\lambda_i^\eps}\right)}{1 + \frac{\vert a_i^\eps -z^\eps\vert}{\lambda_i^\eps}}  \right)
 $$
where $\ds \lim_{R\rightarrow +\infty}\lim_{\eps \rightarrow 0} \delta_{R,\eps} =0$.\\ 

Then, we estimate $I_{ij}^\eps$. Then, there is two cases to consider. First $max(\lambda_i^\eps, \lambda_j^\eps) =o(\vert a_i^\eps -a_j^\eps\vert)$. Then we separate the integral as follow
\be
\begin{split}
I_{ij}^\eps = & \int_{D^+_{ij}} \vert \nabla G(\, .\,,z^\eps) \vert \vert \nabla B_i^\eps \vert \vert \nabla B_{j}^\eps \vert   dz\\
 & + \int_{D^-_{ij}} \vert \nabla G(\, .\,,z^\eps) \vert \vert \nabla B_{i}^\eps \vert \vert \nabla B_{j}^\eps \vert dz,
 \end{split}
\ee
where $D^+_{ij} = \{z \hbox{ s.t. } \langle z -m_{ij}, z_j -z_i \rangle \geq 0 \}$ and $D^-_{ij} = \{z \hbox{ s.t. } \langle z -m_{ij}, z_j -z_i \rangle \leq 0 \}$ and $m_{ij}=\frac{z_i^\eps+ z_j^\eps}{2}$. We are going to estimate the first integral, it would be the same for the second term. But we easily see that 
\be
\int_{D^+_{ij}} \vert \nabla G(\, .\,,z^\eps) \vert \vert \nabla B_i^\eps \vert \vert \nabla B_{j}^\eps \vert   dz \leq t^\eps_{ji} \int_{\R^2} \vert \nabla G(\, .\,,z^\eps) \vert \vert \nabla B_{j}^\eps\vert dz,
\ee
and using  the lemma \ref{estimint1}, we get 
\be
\int_{D^+_{ij}} \vert \nabla G(\, .\,,z^\eps) \vert \vert \nabla B_i^\eps \vert \vert \nabla B_{j}^\eps \vert   dz  =   O\left(t^\eps_{ji} \frac{Ln\left( 2+ \frac{\vert a_j^\eps -z^\eps\vert}{\lambda_j^\eps}\right)}{1 + \frac{\vert a_j^\eps -z^\eps\vert}{\lambda_j^\eps}}  \right) .
\ee
Then we examine the second case, that is to say, up to exchange $i$ and $j$,  $\vert a_i^\eps -a_j^\eps\vert=O(\lambda_j^\eps )$  and $\lambda_i^\eps =o(\lambda_j^\eps)$. Then we easily check that 
 
$$ \vert \nabla B_j^\eps (z)  \vert \leq c\, t^\eps_{ij} \hbox{ on } \R^2 , $$
where $c$ is a positive constant. Hence, we have
\be 
 I_{ij}^\eps \leq O\left(t^\eps_{ij} \int_{\R^2} \vert \nabla G(\, .\,,z^\eps) \vert \vert \nabla B_i^\eps  \vert dz \right). 
 \ee
and using  the lemma \ref{estimint1}, in every case, we have
\beq
\label{Iij2}
 I_{ij}^\eps = O\left(   t^\eps_{ij}  \frac{Ln\left( 2+ \frac{\vert a_i^\eps -z^\eps\vert}{\lambda_i^\eps}\right)}{1 + \frac{\vert a_i^\eps -z^\eps\vert}{\lambda_i^\eps}}  \right) 
+   O\left( t^\eps_{ji} \frac{Ln\left( 2+ \frac{\vert a_j^\eps -z^\eps\vert}{\lambda_j^\eps}\right)}{1 + \frac{\vert a_j^\eps -z^\eps\vert}{\lambda_j^\eps}}  \right)  .
 \eeq
Finally, we estimate $J_i^\eps(z^\eps)$. Firstly we remark that thanks to (\ref{ffa}) we have
\beq
\label{d1}
\vert (u^\eps -B^\eps_i)(z)\vert\leq \sum_{j\not= i}   \vert B^\eps_{j}(z) -B^\eps_j(b_i^\eps)\vert  +  \vert R^\eps(z) -R^\eps(b_i^\eps)\vert.
\eeq
But we easily see that 
 \beq
\label{d2}
\begin{split}
\vert B^\eps_{j}(z) -B^\eps_j(b_i^\eps)\vert &\leq \int_{[a_i^\eps,z]} \vert \nabla B^\eps_{j} \vert dt \\
&=O\left( \int_0^{\vert z-a_i^\eps\vert} \frac{\lambda^\eps_j}{(\lambda^\eps_j)^2 + \vert a_i^\eps + t \tilde{a}_{j}^\eps - a_j^\eps\vert}  dt \right) ,
\end{split}
\eeq
where $\ds \tilde{a}_{i}^\eps =\frac{z-a_i^\eps}{\vert z-a_i^\eps \vert}$. Hence we get 
\beq
\label{d3}
\vert B^\eps_{j}(z) -B^\eps_j(b_i^\eps)\vert \leq t_{ij}^\eps (\vert z-a_i^\eps \vert +\lambda_i^\eps)
\eeq
Moreover we clearly get 
$$   \vert R^\eps(z) -R^\eps(b_i^\eps)\vert \leq r^\eps (\vert z-a_i^\eps \vert +\lambda_i^\eps),$$ 
which gives, with (\ref{d1}) and (\ref{d3}), that
\beq
\label{d33}
\vert (u^\eps -B^\eps_i)(z)\vert\vert \nabla B_i^\eps\vert =O\left( \sum_{j\not= i} t_{ij} ^\eps + r^\eps\right),
\eeq 
and then  
\beq
\begin{split}
\label{d4}
 J_i^\eps(z^\eps) &=O\left(\left(\sum_{j\not= i} t_{ij}^\eps + r^\eps\right) \int_{\R^2} \vert \nabla G(.,z^\eps) \vert \vert \nabla B_i^\eps \vert  dz \right)\\
 &= O\left(\left(\sum_{j\not= i} t_{ij}^\eps + r^\eps\right) \frac{Ln\left( 2+ \frac{\vert a_i^\eps -z^\eps\vert}{\lambda_i^\eps}\right)}{1 + \frac{\vert a_i^\eps -z^\eps\vert}{\lambda_i^\eps}}  \right) .
 \end{split}
\eeq
Then we can estimate $J^\eps$ with respect to the previous terms and$I^\eps$. Indeed
\beq
\label{ef3}
J^\eps(z^\eps)=O\left(\sum_{i=0}^p \frac{\eps^3}{\lambda_i^\eps}\, \frac{Ln\left( 2+ \frac{\vert a_i^\eps -z^\eps\vert}{\lambda_i^\eps}\right)}{1 + \frac{\vert a_i^\eps -z^\eps\vert}{\lambda_i^\eps}}  \right) + \eps^3\sum_{i=0}^p I_i^\eps(z^\eps)  + I^\eps(z^\eps).
\eeq
Hence we just need to estimate $I^\eps(z^\eps)$. In order to do it, we use the weak estimate (\ref{estfaible}) on $\nabla R^\eps$, then we get
\be
I^\eps(z^\eps) = \eps^2 (r^\eps)^\frac{1}{3}\sum_{i=0}^p o\left(\int_{\R^2} \frac{1}{\vert z-z^\eps\vert} \frac{1}{(\lambda^\eps_i+\vert z-a_i^\eps \vert)^\frac{5}{3}}\, dz \right).
\ee
Then we set $y^\eps_i = a^\eps_i -z^\eps$, which gives
\be
I^\eps(z^\eps) = \eps^2 (r^\eps)^\frac{1}{3}\sum_{i=0}^p o\left(\int_{\R^2} \frac{1}{\vert z-y_i^\eps\vert} \frac{1}{(\lambda^\eps_i+\vert z\vert)^\frac{5}{3}}\, dz \right).
\ee
Hence we set the new variable $z=\mu^\eps_i u$ where $\mu^\eps_i= \lambda_i^\eps + \vert y_i^\eps\vert$, then 
\be
I^\eps(z^\eps) = \eps^2 (r^\eps)^\frac{1}{3}\sum_{i=0}^p o\left((\mu^\eps_i)^\frac{-2}{3} \int_{\R^2} \frac{1}{\left\vert u-\frac{y_i^\eps}{\mu^\eps_i} \right\vert} \frac{1}{\left(\frac{\lambda^\eps_i}{\mu^\eps_i}+\vert u\vert\right)^\frac{5}{3}} \, du \right).
\ee
Since the integral are uniformly bounded, we get
\be
I^\eps(z^\eps) = \sum_{i=0}^p o\left(\eps^2 (r^\eps)^\frac{1}{3}(\mu^\eps)^\frac{-2}{3}\right).
\ee
Using the Young inequality
\beq
\label{ZZZ}
\begin{split}
I^\eps(z^\eps) &= \sum_{i=0}^p o\left(  (r^\eps)^\frac{1}{3}  \eps^2(\lambda_i^\eps)^\frac{-2}{3} \left(1 + \frac{\vert y_i^\eps\vert}{\lambda^\eps_i}\right)^\frac{-2}{3}\right)\\
&= \sum_{i=0}^p o\left( r^\eps + \frac{\eps^3}{\lambda_i^\eps}  \left(1 + \frac{\vert y_i^\eps\vert}{\lambda^\eps_i}\right)^{-1} \right) =o\left( r^\eps  + \sum_{i=0}^p\frac{\eps^3}{\lambda^\eps_i}  \frac{1}{1 + \frac{\vert z^\eps -a_i^\eps\vert}{\lambda^\eps_i}}\right).
\end{split}
\eeq

Finally, thanks to (\ref{green2}), (\ref{xxr2}), (\ref{Ij2}), (\ref{Iij2}), (\ref{d4}), (\ref{ef3}) and (\ref{ZZZ}) we get, for any $R>0$, 
\beq
\label{festim}
\begin{split}
\vert \nabla R^\eps (z^\eps)\vert  &= \sum_{i=0}^p  O\left(\left( r_{i,R}^\eps +  \frac{\eps^3}{\lambda_i^\eps} + \sum_{j\not= i} t_{ij}^\eps\right) \frac{Ln\left( 2+ \frac{\vert a_i^\eps -z^\eps\vert}{\lambda_i^\eps}\right)}{1 + \frac{\vert a_i^\eps -z^\eps\vert}{\lambda_i^\eps}}  \right) \\
& +  \delta_{R,\eps} r^\eps.
\end{split}
\eeq
where $\ds \lim_{R\rightarrow +\infty}\lim_{\eps \rightarrow 0} \delta_{R,\eps} =0$.\\ 

Now we are going to show that $r^\eps_{i,R}$ is controlled by $t^\eps$.  First we prove a stronger result when another bubble is closed to the one we consider.\\

{\bf Claim 1: Let $i$ fixed. If there exists $i_0\not= i$ such that $\ds \limsup_{\eps\rightarrow 0} \frac{d_{i_0}^\eps(a_i^\eps)}{\lambda_i^\eps} < + \infty$. Then for every $R>0$, we get 
$$ r^\eps_{i,R} =o( t^\eps) .$$
}

{\it Proof of claim 1 :}\\

Thanks to the previous section, we have
$$\sup_{z\in\R^2} (\min_{0\leq j\leq p} d_j^\eps(z)) \vert \nabla R^\eps \vert =o(1) .$$
Now we fix $R>0$ and $z^\eps \in \Omega_{i,R}^\eps$, we have 
$$ \vert \nabla R^\eps (z^\eps)\vert =o\left( \frac{1}{\min_{j} d_j^\eps(z^\eps)} \right). $$
Then, up to a subsequence, there exists $j_0$ such that
$$ d_{j_0}^\eps(z^\eps) =\min_{0\leq j\leq p} d_j^\eps( z^\eps) ,$$
 and then
 \beq
 \label{rm1} 
 \vert \nabla R^\eps(z^\eps) \vert =o\left( \frac{ d_{j_0}^\eps(z^\eps)}{ (d_{j_0}^\eps(z^\eps)))^2 } \right).
 \eeq
 Moreover, using the fact $z^\eps \in \Omega_{i,R}^\eps$, we easily get
\beq
\label{rm2}
\begin{split}
(\lambda_{i}^\eps )^2 + \vert a_{i}^\eps -a_{j_0}^\eps\vert^2 + (\lambda_{j_0}^\eps)^2 &= O(\vert z^\eps -a_{j_0}^\eps\vert^2 +(\lambda_{j_0}^\eps)^2)\\
& =O ((d_{j_0}^\eps(z^\eps))^2 ).
\end{split}
\eeq
 From the other hand, we get
\beq
\label{rm3}
d_{j_0}^\eps(z^\eps)) \leq d_{i}^\eps(z^\eps)  =O(\lambda_i^\eps)
\eeq
Finally, thanks to (\ref{rm1}), (\ref{rm2}) and (\ref{rm3}), we get
 $$ \vert \nabla R^\eps(z^\eps) \vert =o\left( \frac{\lambda_i^\eps }{(d_{j_0}^\eps(a_{i}^\eps))^2 + (d_{i}^\eps(a_{j_0}^\eps))^2} \right), $$
which proves claim 1.\hfill$\square$\\

Now we prove that if the main estimate is not satisfied, i.e. if $r^\eps$ and $t^\eps$ are not controlled by $\frac{\eps^3}{\lambda^\eps_p}$, then $r^\eps_{i,R}$ is controlled by $t^\eps$. Indeed, else the reminder will be greater than the interactions and will give a non-constant solution of the linearized equation.\\

{\bf Claim 2  :  Either $r^\eps +t^\eps =O\left(\frac{\eps^3}{\lambda_p^\eps}\right)$ or, for any positive number $R$, we have $\ds \max_{0\leq i\leq p} r^\eps_{i,R}  = O(t^\eps)$.}\\  

In particular if $p=0$ we necessary get $t^\eps=0$ and $r^\eps=O\left(\frac{\eps^3}{\lambda_p^\eps}\right)$.\\

{\it Proof of Claim 2 :}\\

Let us assume for contradiction that 
$$ \frac{\eps^3}{\lambda_p^\eps} = o (r^\eps + t^\eps)$$
and that there exists $R$ such that  $\ds t^\eps= o\left(\max_{0\leq i\leq p} r^\eps_{i,R}\right)$. Then, up to a subsequence, we can assume that, there exists $i_0$ such that 
$$r^\eps_{i_0,R}= \max_{0\leq i\leq p} r^\eps_{i,R}. $$
 Of course our hypothesis leads to $t^\eps=o(r^\eps)$. Then we are going to prove that $r^\eps_{i_0,R}=o(r^\eps)$ which, with (\ref{festim}), will give a contradiction and prove the claim.\\

Thanks to the previous claim we can also assume that for every  $j\not= i_0$ we have $\limsup \frac{d_{j}^\eps(a_i^\eps)}{\lambda_{i_0}^\eps} = + \infty$. Then (\ref{cint}) is true and we rescale setting $\tilde{f} = f(\lambda_{i_0}^\eps \, . \, + a_{i_0}^\eps)$.\\

Then, thanks to  (\ref{eqb}), (\ref{eqbc}), (\ref{eqr}), (\ref{eqrc}) and (\ref{cint}), we see that $\tilde{B}_j^\eps$ and $\tilde{R}^\eps$ satisfy the following equations, on every compact subset of $\R^2\setminus\{ S_{i_0}\} $
\be 
\begin{split}
&\vert \nabla \tilde{B}_{i_0}^\eps \vert =O(1), \\
 &\vert \nabla \tilde{B}_{j}^\eps \vert =o(\lambda_{i_0}^\eps r^\eps) \hbox{ for } j\not= i , \\
&\vert \nabla \tilde{R}^\eps \vert =O \left( \lambda_{i_0}^\eps r^\eps \right),
\end{split}
\ee
and
\be
\Delta\left( \sum_{j=0}^k \tilde{B}_j^\eps +  \tilde{R}^\eps \right)= -2 \left( \sum_{j=0}^k \tilde{B}_j^\eps +  \tilde{R}^\eps \right)_x \wedge \left( \sum_{j=0}^k \tilde{B}_j^\eps +  \tilde{R}^\eps \right)_y + o\left(\lambda_{i_0}^\eps r^\eps \right) ,
\ee
and  the relation of quasi-conformality
\be
\begin{split}
&\left\langle \sum_{j=0}^k  (\tilde{B}_j^\eps)_x + \tilde{R}^\eps_x , \sum_{j=0}^k (\tilde{B}_j^\eps)_y +  \tilde{R}^\eps_y \right\rangle = o\left( \lambda_{i_0}^\eps r^\eps \right) ,\\
&\left\langle \sum_{j=0}^k  (\tilde{B}_j^\eps)_x +  \tilde{R}^\eps_x , \sum_{j=0}^k (\tilde{B}_j^\eps)_x +  \tilde{R}^\eps_x \right\rangle - \left\langle \sum_{j=0}^k  (\tilde{B}_j^\eps)_y +  \tilde{R}^\eps_y , \sum_{j=0}^k (\tilde{B}_j^\eps)_y +  \tilde{R}^\eps_y \right\rangle = o \left(\lambda_{i_0}^\eps r^\eps \right),
\end{split}
\ee
and the initial conditions
\be
\begin{split} 
& \nabla\left( \sum_{j\not= i_0} \tilde{B}_j^\eps +  \tilde{R}^\eps \right)(0)=o  ( \lambda_{i_0}^\eps r^\eps ),\\
& \nabla^2\left( \sum_{j\not= i_0} \tilde{B}_j^\eps +  \tilde{R}^\eps \right)(\nabla \tilde{B}_{i_0}^\eps)(0)=o( \lambda_{i_0}^\eps r^\eps ).
\end{split}
\ee
Then, thanks to standard elliptic theory, see \cite{GT}, we get that   $\frac{\tilde{R}^\eps}{\lambda_{i_0}^\eps r^\eps}$ converge in $C^2_{loc} (\R^2\setminus\{ S_{i_0} \} )$ to  $\tilde{R}$ which satisfies
\be
\Delta \tilde{R}= -2 \left( \omega^{i_0}_x \wedge \tilde{R}_y  +  \tilde{R}_x \wedge \omega^{i_0}_y\right) ,
\ee
and the relations of conformality ,
\be
\begin{split}
&\langle \omega^{i_0}_x ,  \tilde{R}_y \rangle + \langle  \tilde{R}_x , \omega^{i_0}_y  \rangle = 0 ,\\
&\langle \omega^{i_0}_x ,  \tilde{R}_x \rangle - \langle \omega^{i_0}_y , \tilde{R}_y  \rangle = 0
\end{split}
\ee
 and 
 \be
\begin{split}
&\nabla \tilde{R}(0)=0 ,\\
&\nabla^2 \tilde{R}(\nabla \omega^{i_0})(0)=0 .
\end{split}
\ee
Moreover, $\nabla \tilde{R}$ is uniformly bounded on $\R^2\setminus \{ S_{i_0} \}$, then it can be extended to a smooth function of $\R^2$ which satisfies the same equation and whose gradient is still uniformly bounded. Finally applying proposition \ref{lin}, we see that 
 $$ \nabla \tilde{R} \equiv 0,$$ 
which proves that $r^\eps_{i_0,R}=o(r^\eps)$. As already said this last estimate contradicts (\ref{festim}) applied with $R$ big enough, $\eps$ small enough and $z^\eps$ such that $\nabla R^\eps (z^\eps)=\frac{r^\eps}{2}$, which finally proves claim 2. \hfill $\square$\\

Finally applying the fundamental estimate (\ref{festim}) with $z^\eps$ such that $\nabla R^\eps (z^\eps)=\frac{r^\eps}{2}$, $R$ big enough and $\eps$ small enought, we get thanks to claim  2, that
\beq
\label{festim2}
r^\eps  =  O\left(t^\eps + \frac{\eps^3}{\lambda_p^\eps} \right),
\eeq
In order to get our desired estimate, it suffices to prove that $t^\eps =O\left( \frac{\eps^3}{\lambda_p^\eps}\right)$. This fact is postponed to the last section. Now we will take it  as  proved while proving the theorem. We can remark that this estimate is automatically satisfies when there is only one bubble, since the interaction term vanishes. 

%%%%%%%%%%%%%%%%%%%%%%%%%%%%%%%%%%%%%%%%%%%%%%%%%
\section{Proof of theorem \ref{thmain} }
\label{s4}
In this section we assume that 
$$r^\eps + t^\eps =O\left(\frac{\eps^3}{\lambda_p^\eps} \right). $$
We are going to use this estimate looking at the highest bubble, that is to say $\omega^p$.\\

 We set $\tilde{f}=f(\lambda_p^\eps\; . + a_p^\eps)$, thanks to (\ref{eqr}), $\tilde{R}^\eps = \tilde{u}^\eps -\tilde{B}_p^\eps$ then satisfies, on every compact set of $\R^2$, 
\beq
\label{last}
\begin{split}
\Delta \tilde{R}^\eps &= -2 \sum_{i=0}^{p-1}(\tilde{B}_p^\eps)_x \wedge \tilde{R}^\eps_y +  \tilde{R}^\eps_x \wedge (\tilde{B}_p^\eps )_y \\
&+\eps^3 \left( \frac{1}{6} \mathrm{Ric}_{ij,k}(\tilde{B}_p^\eps)^i(\tilde{B}_p^\eps)^j(\tilde{B}_p^\eps)^k ((\tilde{B}_p^\eps)_x \wedge (\tilde{B}_p^\eps) _y)  \right. \\
&+\frac{1}{3} R_{ikmj,n}(\tilde{B}_p^\eps)^k(\tilde{B}_p^\eps)^m(\tilde{B}_p^\eps)^n ((\tilde{B}_p^\eps)_x \wedge(\tilde{B}_p^\eps)_y )^i \\
&+\left.  B_{ijkmn}(\tilde{B}_p^\eps)^m(\tilde{B}_p^\eps)^n \left\langle\nabla (\tilde{B}_p^\eps)^{i},\nabla (\tilde{B}_p^\eps)^{k}\right\rangle \right)\\
&  +\tilde{R}^\eps_x \wedge \tilde{R}^\eps_y  + o\left( \sum_{i=0}^p \vert \nabla \tilde{B}_i^\eps \vert  \vert \nabla \tilde{R}^\eps \vert  \right) + O( \eps^4 \vert \nabla \tilde{u}^\eps\vert^2)
\end{split}
\eeq
where $B_{ijkmn}$ is defined in appendix \ref{a2}. Then, dividing (\ref{last}) by $\eps^3$ and thanks to the standard elliptic theory, see \cite{GT}, up to a subsequence, $\frac{\tilde{R}^\eps}{\eps^3}$ converge in $C^2_{loc}(\R^2)$ to $\tilde{R}$ solution on $\R^2$ of
\beq
\label{fin}
\begin{split}
\Delta \tilde{R} &= -2 (\hat{\omega}^p_x \wedge  \tilde{R}_y  +  \tilde{R}_x \wedge \hat{\omega}^p _y ) +\frac{1}{6} \mathrm{Ric}_{mn,k}(\hat{\omega}^p + p_p)^m(\hat{\omega}^p + p_p)^n(\hat{\omega}^p + p_p)^k (\hat{\omega}^p_x \wedge \hat{\omega}^p _y )_j \\
&+\frac{1}{3} R_{ikmj,n}(\hat{\omega}^p + p_p)^k(\hat{\omega}^p + p_p)^m(\hat{\omega}^p + p_p)^n (\hat{\omega}^p_x \wedge \hat{\omega}^p _y )^i \\
&+  B_{ijkmn}(\hat{\omega}^p + p_p)^m(\hat{\omega}^p+p_p)^n \left\langle\nabla (\hat{\omega}^p)^{i},\nabla (\hat{\omega}^p)^{k}\right\rangle.
\end{split}
\eeq  
Then, up to compose  $\tilde{R}$ with an homography, we can assume that  $\hat{\omega^p } =\omega$. We also replace $p_p$ by $p$.\\ 

Moreover, we know that $\omega_x$, $\omega_y$ and $x\, \omega_x + y\, \omega_y$ are solution of the linearized operator, then testing (\ref{fin}) against these functions we should find some informations.  From now to the end of the proof we denote $\omega_x$, $\omega_y$ and $x\, \omega_x + y\, \omega_y$ by $Y^1$,$Y^2$ and $Y^3$. Let $R>0$, then we have
\be
\begin{split}
\int_{B(0,R)} Y^l \Delta \tilde{R}\, dz& =\int_{B(0,R)}   -2 \langle Y^l,\omega_x \wedge  \tilde{R}_y  +  \tilde{R}_x \wedge \omega _y \rangle\, +C_{j}(p_\infty,z) (Y^l)^j\, dz ,
\end{split}
\ee
where 
\be 
\begin{split}
C_{j}(p_\infty,z) &= B_{ijkmn}(p_\infty)(\omega+p)^m(\omega+p)^n \left\langle\nabla \omega^{i},\nabla \omega^{k}\right\rangle \\
&+\frac{1}{3} R_{ikmj,n}(p_\infty)(\omega + p)^k(\omega + p)^m(\omega + p)^n (\omega_x \wedge\omega_y )^i\\
& -\frac{1}{6} \mathrm{Ric}_{ij,k}(p_\infty)(\omega+p)^i(\omega+p)^j(\omega+p)^k (\omega_x \wedge\omega_y ) .
 \end{split}
 \ee
Integrating by parts, we get 
\beq
\label{lint}
\begin{split}
&\int_{B(0,R)}  \langle \tilde{R} , \Delta Y^l +2 (Y^l_{x} \wedge  \omega_{y}  +  \omega_{x} \wedge Y^l_y )\rangle dz = \int_{B(0,R)}  C_{j}(p_\infty,z)(Y^l)^j dz\\
&+ O\left( \int_{\partial B(0,R)}  (\vert \nabla \omega\vert \vert Y^l\vert + \vert \omega\vert \vert\nabla  Y^l\vert +  \vert\tilde{R}\vert \vert \nabla  Y^l\vert +  \vert \nabla \tilde{R}\vert   \vert Y^l\vert)dz\right) .
\end{split}
\eeq
Thanks to the fact that $\nabla \tilde{R}$ satisfies (\ref{festim}), we get that
\be
\begin{array}{cc} 
\begin{cases}
\vert \tilde{R} \vert =o(\vert z\vert)  \\
  \vert \nabla \tilde{R} \vert =o \left(1\right)
  \end{cases} 
&\hbox{ as } z\rightarrow + \infty .  \end{array}
\ee
Moreover, thanks to the formulas of section 2, we also get the following estimates
\be
\begin{array}{cc} 
\begin{cases}
 \vert \omega\vert =O(1) \\
  \vert \nabla \omega\vert =O \left(\frac{1}{\vert z\vert^2 }\right)
  \end{cases} 
&\hbox{ as } z\rightarrow + \infty ,  \end{array}
\ee
and 
\be
\begin{array}{cc} 
\begin{cases}
 \vert Y^k\vert =O \left(\frac{1}{\vert z\vert }\right)  \\
  \vert \nabla Y^k \vert =O \left(\frac{1}{\vert z\vert^2 }\right) 
  \end{cases}
&\hbox{ as } z\rightarrow + \infty .  \end{array}
\ee
Thanks to these estimates, passing to the limit in (\ref{lint}) as $R$ goes to infinity, we get 
\be
\begin{split}
 &\int_{\R^2}  B_{ijlm}(p_\infty)(\omega+p)^m(\omega+p)^n \left\langle\nabla \omega^{i},\nabla \omega^{k}\right\rangle(Y^l)^j  dz  =\\
 &\frac{1}{6} \int_{\R^2} \mathrm{Ric}_{ij,k}(p_\infty)(\omega+p)^i(\omega+p)^j(\omega+p)^k (\omega_x \wedge\omega_y )(Y^l)^j dz \\ 
 &- \frac{1}{3} \int_{\R^2} R_{ikmj,n}(p_\infty)(\omega + p)^k(\omega + p)^m(\omega + p)^n (\omega_x \wedge\omega_y )^i (Y^l)^j  dz \,.
\end{split}
\ee
Then changing the variable via $y=\omega(z)$ and  using (\ref{nablai}), we get the following integral on the sphere
\be
\begin{split}
 &\int_{\S^2}  B_{ijlm}(p_\infty)(y+p)^m(y+p)^n (\delta^{ik}-y^i y^k)(Y^l)^j  dv_h  =\\
 &\frac{1}{6}  \int_{S^2}\mathrm{Ric}_{ij,k}(p_\infty)(y+p)^i(y+p)^j(y+p)^k y_j(Y^l)^j dv_h \\ 
 &-\frac{1}{3} \int_{S^2} R_{ikmj,n}(p_\infty)(y + p)^k(y + p)^m(y + p)^n y^i (Y^l)^j  dv_h\, .
\end{split}
\ee
Then we compute $Y^l(y)$, thanks to the formulas of section 2, we get that 
\be
\begin{split}
&\omega_x(\pi(y)) = \left(\begin{array}{c}1 \\0 \\0\end{array}\right) +\left(\begin{array}{c} -y^3 \\0 \\ y^1\end{array}\right) -y^1 y , \\
&\omega_y(\pi(y)) = \left(\begin{array}{c}0 \\1 \\0\end{array}\right)+ \left(\begin{array}{c} 0 \\ -y^3 \\ y^2\end{array}\right) - y^2 y, \\
&(x\, \omega_x + y\, \omega_y)(\pi(y)) = \left(\begin{array}{c}0 \\0 \\1\end{array}\right) - y^3 y,
\end{split}
\ee
Now taking in account that every integrand with an odd number of $y$ vanishes and using the symmetry of the Riemannian tensor, we get, after a straitforward computation, that  
$$ \int_{S^2} B_{ijkmn}(p_\infty) y^m y^n(\delta^{ik}-y^{i} y^{k})(\delta^{jl} -y^j y^l) dv_h=0  \; \hbox{ for all } l.$$
We remark that this expression is independent of $p$. It is naturual since $p$ depends on the center of chart we have choose at the beginning of our analysis. But this center of chart has been chosen arbitrary, hence the result mustn't depend on $p$.\\

Finally, replacing $B_{ijkmn}$ by its expression, we get that
\be
\int_{S^2} ( 4R_{kmij,n}(p_\infty) + 2 R_{imnj,k}(p_\infty) - R_{imnk,j}(p_\infty) ) y^m y^n(\delta^{ik}-y^{i} y^{k})(\delta^{jl} -y^j y^l) dv_h= 0
\ee
We have the following standard formulas on the sphere
$$ \int_{S^2} y^m y^n dv_h =\frac{4\pi}{3} \delta^{mn} \hbox{ and }   \int_{S^2} y^m y^n y^l y^j dv_h =\frac{4\pi}{15}(\delta^{mn}\delta^{jl} + \delta^{mj}\delta^{nl} +\delta^{ml}\delta^{nj}) $$
which gives that 
\be
 \mathrm{Ric}_{ml}^{\ \ \ ,m}(p_\infty)  =0   \hbox{ for all }l  .
\ee
Finally, thanks to the second Bianchi identity, we have
$$ \nabla \mathrm{Scal} (p_\infty)=0 .$$
This  achieves the proof of  the theorem.\hfill$\blacksquare$
%%%%%%%%%%%%%%%%%%%%%%%%%%%%%%%

\section{estimate on the bubble interaction}
\label{s7}
The aim of this section is to prove the following claim.\\

{\bf Claim : $t^\eps=O(\frac{\eps^3}{\lambda^\eps_p})$.}\\

{\it Proof of the claim :}\\

We assume for contradiction that $\frac{\eps^3}{\lambda^\eps_p}=o(t^\eps)$. First we remark that this implies thanks to (\ref{festim}) that 
$$r^\eps =O(t^\eps).$$

Before we start the proof , we give some complementary definition on the interaction.\\

Let $I=\{ i \vert \exists j \hbox{ s.t. }  \liminf \frac{ t^\eps_{ij}}{ t^\eps} >0 \}$ be the set of indices whose bubbles receive a maximal interaction and $T_i =\{ j \hbox{ s.t. } \liminf \frac{ t^\eps_{ij}}{ t^\eps} >0  \}$ be the set of indices whose bubbles give this maximal interaction.\\

First we prove that  each element of $I$ received at least two maximal interaction.\\ 

{\bf Claim 1 : For all $i_0 \in I $  we have
\beq
\label{rem}
\vert T_{i_0}\vert >1.
\eeq
}\\
{\it Proof of Claim 1 :}\\

Let us assume for contradiction that  $T_{i_0}=\{j_0\}$. Then, we prove that  if a bubble "contains " $B_{i_0}^\eps$, it can't receive any maximal interaction. Indeed, else $B_{i_0}^\eps$ would received more than one maximal interaction which contradict our hypothesis.\\

{\bf Claim 1.1 : Let $i\not= i_0$ then either $\ds \limsup_{\eps\rightarrow 0} \frac{d_{i_0}^\eps(a_{i}^\eps)}{\lambda_{i}^\eps} =+\infty $ or $t_{ik}^\eps =o(t^\eps)$ for all $k\not= i$.}\\

{\it Proof of Claim 1.1 :}\\

Let $i\not= i_0$ and let us assume that  $\limsup_{\eps\rightarrow 0} \frac{d_{i_0}^\eps(a_{i}^\eps)}{\lambda_{i}^\eps} <+\infty $. Then thanks to (\ref{orthof}), we have
\beq
\label{rr2}
\begin{cases}
\lambda_{i_0}^\eps = o(\lambda_i^\eps) ,\\
\vert a_i^\eps - a_{i_0}^\eps \vert \leq \lambda_i^\eps .
\end{cases}
\eeq
Let us assume for contradiction assume that there exists $k\not=  i $ such that $t^\eps = O(t_{ik}^\eps)$. Remarking  that  we necessarily get  that  $t^\eps_{ki} = O(t_{ik}^\eps)$, then we have
\beq
\label{rr1}
\limsup_{\eps\rightarrow 0} \frac{d_{k}^\eps(a_i^\eps)}{\lambda_i^\eps} =+\infty .
\eeq
Else we have
\be
\begin{cases}
\lambda_{k}^\eps = o(\lambda_i^\eps) ,\\
\vert a_i^\eps - a_{k}^\eps\vert =O(  \lambda_i^\eps) ,
\end{cases}
\ee
which gives that
\be
\begin{cases}
t_{ik}^\eps = O\left( \frac{\lambda_k^\eps}{(\lambda_i^\eps)^2} \right),\\
 \frac{1}{\lambda_i^\eps} = O\left( t_{ki}^\eps\right) ,
\end{cases}
\ee
and leads to $t^\eps_{ik} = o(t_{ki}^\eps)$, which is clearly a contradiction and proves (\ref{rr1}).\\
Then, thanks to  (\ref{rr2}) and  (\ref{rr1}), we also get
\be
\begin{split}
&t_{ik}^\eps = \frac{\lambda_k^\eps}{(\lambda_k^\eps)^2 + (\lambda_{i}^\eps)^2  +  \vert a_{i}^\eps - a_{k}^\eps \vert^2 )}=  O\left(\frac{\lambda_k^\eps}{(d_k^\eps (a_i^\eps))^2}   \right) =O\left(\frac{1}{d_k^\eps (a_i^\eps)}   \right)
,\\
&\hbox{ and } \\
&\frac{1}{\lambda_i^\eps} = O\left(\frac{\lambda_i^\eps}{(\lambda_{i_0}^\eps)^2 + (\lambda_{i}^\eps)^2  +  \vert a_{i_0}^\eps - a_{i}^\eps \vert^2 } \right)= O(t_{i_0i}^\eps),
\end{split}
\ee
Then, thanks to (\ref{rr1}), we easily get that $t^\eps = o(t_{i_0 i}^\eps)$, which is a contradiction and achieves the proof.\hfill$\square$\\

Now we are going to give a decreasing estimate on $\nabla R^\eps$ around $a_{i_0}^\eps$. Let $R>0$ and $z^\eps$ such that $\vert z^\eps -a_{i_0}^\eps \vert = R\lambda_{i_0}^\eps$. Thanks to  (\ref{festim})  we have
\beq
\label{ffestim}
\vert \nabla R^\eps (z^\eps)\vert  \leq o(t^\eps) + \left(\sum_{i=0}^p \sum_{j\not= i} t_{ij}^\eps\right) O\left(\frac{Ln\left( 2+ \frac{\vert a_i^\eps -z^\eps\vert}{\lambda_i^\eps}\right)}{1 + \frac{\vert a_i^\eps -z^\eps\vert}{\lambda_i^\eps}}  \right),
\eeq
Then there is two cases to consider. Let $i\not= i_0$, if $\ds \limsup_{\eps\rightarrow 0} \frac{d_{i_0}^\eps(a_i^\eps)}{\lambda_{i}^\eps} =+\infty$ then we easily see that 
\beq
\label{last1}
\frac{Ln\left( 2+ \frac{\vert a_i^\eps -z^\eps\vert}{\lambda_i^\eps}\right)}{1 + \frac{\vert a_i^\eps -z^\eps\vert}{\lambda_i^\eps}}  =o(1).
\eeq
 Else $\ds \limsup_{\eps\rightarrow 0} \frac{d_{i_0}^\eps(a_i^\eps)}{\lambda_{i}^\eps} <+\infty$ and thanks to claim 1.1,  for all $j\not= i_0$, we get that
 \beq
 \label{last2}
 t^\eps_{ij} =o(t^\eps) .
 \eeq
 Finally thanks to (\ref{last1}) and (\ref{last2}) we get 
\be
\vert \nabla R^\eps (z^\eps)\vert  \leq o(t^\eps) + O\left(t^\eps \frac{Ln\left( 2+R\right)}{1 +R}  \right).
\ee
So we see that $\frac{\nabla R^\eps}{t^\eps}$ decreases at infinity, so it cannot compensate $\frac{\nabla B_{j_0}^\eps}{t^\eps}$ which is constant. But the sum of its two function should goes to a solution of the linearized equation, that is to say zero which will leads us to a contradiction. The following is devoted to prove what we have just sketched.\\

First we prove, for all $j\not=i_0$, that 
\beq
\label{hypc1}
\lim_{\eps\rightarrow 0} \frac{d_j^\eps(a_{i_0}^\eps)}{\lambda_{i_0}^\eps} \not= 0 .
\eeq
Indeed, thanks to (\ref{orthof}), we get that 
\be
\lambda_{i_0}^\eps \lambda_{j_0}^\eps = o ((d_{i_0}^\eps(a_{j_0}^\eps))^2 + (d_{i_0}^\eps(a_{j_0}^\eps))^2) ,
\ee
if (\ref{hypc1}) doesn't hold for some $j$, then  we get that
$$ t_{i_0 j_0}^\eps =o\left(\frac{1}{\lambda_{i_0}^\eps}\right) =o(t_{ji_0}^\eps),$$
which is a contradiction and proves  (\ref{hypc1}). Hence it  insure that (\ref{ladjust}) is satisfied.\\

Then we rescale the bubble $i_0$ setting $\tilde{f}= f(\lambda_{ i_0}^\eps z+a_{i_0}^\eps)-f(a_{i_0}^\eps)$.  Thanks to our assumption (\ref{eqb}), (\ref{eqbc}), (\ref{eqr}) and  (\ref{eqrc}),  we see that $\tilde{B}_i^\eps$ and $\tilde{R}^\eps$ satisfy the following equations, on every compact subset of $\R^2\setminus\{ S_{i_0}\} $,
\be 
\begin{split}
&\vert \nabla \tilde{B}_{i_0}^\eps \vert =O(1) ,\\
&\vert \nabla \tilde{B}_{j_0}^\eps \vert =O(\lambda_{i_0}^\eps t^\eps) ,\\
& \vert \nabla \tilde{B}_{i}^\eps \vert =o(\lambda_{i_0}^\eps t^\eps) \hbox{ for } i\not\in \{ i_0, j_0\}, \\
&\vert \nabla \tilde{R}^\eps \vert =O \left( \lambda_{i_0}^\eps t^\eps \right) .
\end{split}
\ee
Then, we also get
\be
\begin{split}
\Delta \tilde{B}_i^\eps &= -2(\tilde{B}_i^\eps)_x \wedge (\tilde{B}_i^\eps)_y + O(\eps^3) \hbox{ for all } i,\\
\Delta(\tilde{B}_{j_0}^\eps +  \tilde{R}^\eps )&= -2 ( (\tilde{B}_{i_0}^\eps)_x \wedge (\tilde{B}_{j_0}^\eps +  \tilde{R}^\eps )_y + \wedge (\tilde{B}_{j_0}^\eps +  \tilde{R}^\eps )_x \wedge (\tilde{B}_{i_0}^\eps)_y) + o(\lambda_{i_0}^\eps t^\eps ),
\end{split}
\ee
 and the relation of quasi-conformality
\be
\begin{split}
&\langle (\tilde{B}_{i_0}^\eps )_x , (\tilde{B}_{j_0}^\eps + \tilde{R}^\eps )_y \rangle+ \langle (\tilde{B}_{j_0}^\eps + \tilde{R}^\eps )_x,  (\tilde{B}_{i_0}^\eps )_y \rangle =  o(\lambda_{i_0}^\eps t^\eps ), \\
& \langle(\tilde{B}_{i_0}^\eps )_x  , (\tilde{B}_{j_0}^\eps + \tilde{R}^\eps )_x \rangle- \langle(\tilde{B}_{j_0}^\eps+  \tilde{R}^\eps )_y ,(\tilde{B}_{i_0}^\eps )_y \rangle =  o(\lambda_{i_0}^\eps t^\eps ).
\end{split}
\ee
Moreover, thanks to (\ref{ladjust}), we have
 \beq
 \label{init}
\begin{split}
&\nabla (\tilde{B}_{j_0}^\eps + \tilde{R}^\eps )(0)=o(\lambda_{i_0}^\eps t^\eps) ,\\
&\nabla^2 (\tilde{B}_{j_0}^\eps + \tilde{R}^\eps) (\nabla \tilde{B}_{i_0}^\eps ) =o(\lambda_{i_0}^\eps t^\eps) .
\end{split}
\eeq
Then, thanks to standard elliptic theory, see \cite{GT}, we get that   $\frac{\tilde{B}^\eps_{j_0} +\tilde{R}^\eps}{\lambda_{i_0}^\eps t^\eps}$  converges in $C^2_{loc} (\R^2\setminus\{ S_{i_0} \} )$ to $\tilde{S}$ which satisfies
\be
\Delta \tilde{S}= -2 \left( \omega^{i_0}_x \wedge \tilde{S}_y  +  \tilde{S}_x \wedge \omega^{i_0}_y\right),
\ee
with  the relations of conformality
\be
\begin{split}
&\langle \omega^{i_0}_x , \tilde{S}_y \rangle + \langle  \tilde{S}_x, \omega^{i_0}_y \rangle  = 0 ,\\
& \langle \omega^{i_0}_x , \tilde{S}_x \rangle- \langle \tilde{S}_y , \omega^{i_0}_y \rangle = 0,
\end{split}
\ee
 and the initial data,
 \be
\begin{split}
&\nabla\tilde{S}(0)=0 ,\\
&\nabla^2 \tilde{S} (\nabla \omega^{i_0}(0)) =0 .
\end{split}
\ee
From the other hand, using the fact that $t^\eps_{j_0 i_0} = O(t^\eps_{i_0 j _0})$  and (\ref{orthof}), we see that
$$\lim_{\eps\rightarrow 0} \frac{d_{j_0}(a_{i_0}^\eps)}{\lambda_{i_0}} =+\infty .$$
Then, we  deduce that  $\frac{\nabla \tilde{B}_{j_0}^\eps}{\lambda_{i_0}^\eps t^\eps}$ is  uniformly bounded and satisfies
$$ \Delta \tilde{B}_{j_0}^\eps =o((\lambda_{i_0}^\eps t^\eps)^2)$$
Hence, we easily deduce that $\frac{\nabla \tilde{B}_{j_0}^\eps}{\lambda_{i_0}^\eps t^\eps}$ converges to a constant vector different from zero on $\R^2$. Moreover $\nabla \tilde{R}$ is uniformly bounded on $\R^2\setminus \{ S_{i_0} \}$, then it can be extend to a smooth function  of $\R^2$ which satisfies the same equation and whose gradient is still uniformly bounded. Now we can apply  lemma \ref{lin} to $\tilde{S} = \tilde{R} +  \tilde{B}_{j_0}$ and  we get that  
 $$ \nabla (\tilde{R} +  \tilde{B}_{j_0})\equiv 0,$$ 
 which proves that, for $R$ big enought and $\eps$ small enough, we have
 $$ \vert \nabla R^\eps (z^\eps)\vert \geq \frac{t^\eps}{2}.$$
Then we have a contradiction with (\ref{ffestim}) which achieves the proof of claim 1.\hfill$\square$\\

Now our aim is to find among all bubbles with a maximal interaction a good configuration, that is to say one  where the bubbles are separated. Then passing  to the  limit we will get a contradiction. Indeed,  we will get a sum of plane which will be minimal which can't be approximate by embedded surfaces. This planes are what is seen of $B^\eps_j $ from $a_i^\eps$ asymptoticly, that is to say a tangent plane.\\ 

{\bf Claim 2 : There exits $i_0\in I$ such that, setting $d^\eps = \min \{ d_{j}^\eps(a_{i_0}^\eps) \hbox{ s.t. } j \in T_{i _0}\}$,  for  $k\in T_{i_0}$ either $\ds \lim_{\eps\rightarrow 0 } \frac{\lambda_{k}^\eps}{ \lambda_{i_0}^\eps} >0$ or  $d^\eps =o(d_{i_0}^\eps(a_k^\eps))$.}\\

 That is to say either a bubble of $T_{i_0}$ is  at $d_\eps$ from $ a_{i_0}^\eps$ and have a reciprocal  interaction with the bubble  $B_{i_0}^\eps$ or goes to infinity.\\ 

{\it Proof of claim 2 : }\\

We are going to find $i_0\in I$ which satisfies our claim by induction. In fact let $i_0 \in I$, there is two possibilities :\\

{\bf First case},  there exists $j_0\in T_{i_0}$ such that  $\ds \lim_{\eps\rightarrow 0} \frac{t_{j_0i_0}^\eps}{t^\eps}>0$, that is to  say the interaction between $i_0$ and $j_0$ is reciprocal, and we get 
\beq
\label{equi}
\lim_{\eps\rightarrow 0} \frac{\lambda_{i_0}^\eps}{\lambda_{j_0}^\eps} >0.
\eeq
Then, for every $k \in T_{i_0}$, either 
 $$\lim_{\eps\rightarrow 0}  \frac{d_{k}^\eps(a_{i_0}^\eps)}{d_{j_0}^\eps(a_{i_0}^\eps)}=+\infty,$$ 
 or, thanks to (\ref{equi}), we get 
\be
\begin{split}
t_{i_0 j_0} &=O\left( \frac{\lambda_{i_0}^\eps}{(d_{j_0}^\eps(a_{i_0}^\eps))^2}\right)=O\left( \frac{\lambda_{i_0}^\eps}{(d_{k}^\eps(a_{i_0}^\eps))^2 + (d_{j_0}^\eps(a_{i_0}^\eps))^2} \right) \\
&= O\left( \frac{\lambda_{i_0}^\eps}{(d_{i_0}^\eps(a_{k}^\eps))^2 + (d_{k}^\eps(a_{i_0}^\eps))^2 } \right),
\end{split}
\ee
which gives
 $$\lim_{\eps\rightarrow 0} \frac{t_{k i_0}^\eps}{t^\eps}>0 ,$$
and finally we get  
$$ \lim_{\eps\rightarrow 0 } \frac{\lambda_{k}^\eps}{ \lambda_{i_0}^\eps} >0.$$
Hence, in that cases $i_0$ satisfies our claim.\\

{\bf Second Case}, else, for all $j\in T_{i_0}$, we have  $\ds \lim_{\eps\rightarrow 0} \frac{t_{ji_{0}}^\eps}{t^\eps}=0$. Let  $j_0 \in T_{i_0}$ such that  $\lambda_{j_0}^\eps \leq \lambda_{k}^\eps $ for all $k \in T_{i_0} $. Then  $j_0$ satisfies  
\be
1<\vert T_{j_0} \vert < \vert T_{i_0} \vert.
\ee
In fact, let us prove that  $0<\vert T_{j_0} \vert  $ then the right hand side inequality will follow from claim 1. Let  $k\in T_{i_0}\setminus \{ j_0 \}$, which is not empty thanks to claim 1. Thanks to the fact that $t_{k i_{0} }^\eps =o(t^\eps)$, we have 
\beq
\label{r11}
\lambda_{i_{0} }^\eps =o(\lambda_k^\eps).
\eeq
Moreover, thanks to our hypothesis on $j_0$ and $k$, that is to say 
\be 
\lambda_{j_0}^\eps \leq \lambda_{k}^\eps \hbox{ and } \lim_{\eps\rightarrow 0 }\frac{t_{i_{0} k}^\eps}{t^\eps_{i_0 j_0}} >0,
\ee
we have
\beq
\label{r22}
d_{j_0}^\eps(a_{i_0}^\eps)= O(d_{k}^\eps(a_{i_0}^\eps)).
\eeq
Then, thanks to (\ref{r11}) and (\ref{r22}), we get 
\be
\begin{split}
 t_{ i_0 k}^\eps &=O\left(\frac{\lambda_k^\eps}{(d_{k}^\eps(a_{i_0}^\eps))^2} \right) \\
&=O\left(\frac{\lambda_k^\eps}{(d_{k}^\eps(a_{i_0}^\eps))^2 +(\lambda_{j_0}^\eps)^2} \right) \\
&=O\left(\frac{\lambda_k^\eps}{(d_{k}^\eps(a_{j_0}^\eps))^2 +(d_{j_0}^\eps(a_k^\eps))^2} \right), 
\end{split}
\ee
which proved that $t_{ i_0k}^\eps = O(t_{ j_0k}^\eps )$ and  the left hand side of the desired inequality.\\

In order to show the right hand side inequality, we show that $ T_{j_0}\subset T_{i_0}\setminus \{j_0\}$ . Indeed let $k$ in the complementary of $T_{i_0}\setminus \{j_0\}$, then $d_{j_0}^\eps(a_{i_0}^\eps) =O(d_{j_0}^\eps (a_{k}^\eps))$, else  using (\ref{r11}) we easily get that $t_{i_0 j_0}^\eps=o(t_{kj_0}^\eps)$ which is absurd. Hence we have 
\be
\begin{split}
t_{j_0k}^\eps &= O\left(\frac{\lambda_k^\eps}{(\lambda_k^\eps)^2 + (\lambda_{j_0}^\eps)^2  + \vert  a_{j_0}^\eps - a_k^\eps  \vert^2}  \right)\\
&= O\left(\frac{\lambda_k^\eps}{(\lambda_k^\eps)^2 + (\lambda_{i_0}^\eps)^2 + \vert a_{i_0}^\eps -a_{j_0}^\eps \vert^2 + \vert  a_{j_0}^\eps - a_k^\eps  \vert^2}  \right)\\
&= O\left(\frac{\lambda_k^\eps}{(\lambda_k^\eps)^2 + (\lambda_{i_0}^\eps)^2 + \vert a_k^\eps -a_{i_0}^\eps \vert^2}  \right)\\
 &=O(t_{i_0 k}^\eps) \\
 &=o(t^\eps) ,
 \end{split}
 \ee 
which proves the assertion.\\

Hence if $i_0$ doesn't satisfy the claim we restart with $j_0$, then this induction achieve since the sequence $\vert T_{i_0}\vert $ is strictly decreasing and greater than $1$, which proves the claim.\hfill$\square$ \\

Now, we are in position to prove the main claim of this section. Let $i_0$ as in the previous claim, then  we rescale the space around $a_{i_0}^\eps$ setting $\tilde{f}= \frac{f(d^\eps z+a_{i_0}^\eps)}{d^\eps t^\eps}$, where $d^\eps =\min\{ \vert a_{i_0}^\eps -a_ j^\eps \vert \hbox{ s.t. } j\in T_{i_0} \}$. Thanks to (\ref{eqb}), (\ref{eqbc}), (\ref{eqr}), (\ref{eqrc})  then  $\tilde{B}_i^\eps$ and $\tilde{R}^\eps$ satisfy the following equations, on every compact subset of $\R^2 \setminus \tilde{S}_{i_0}$, where $\ds \tilde{S}_{i_0}=\lim_{\eps\rightarrow 0} \left\{ \frac{ a_{j}^\eps -a_{i_0}^\eps}{d^\eps} \hbox{ s.t. } 1\leq j\leq k \right\}$, 
\be 
\begin{split}
&\vert \nabla \tilde{B}_{i}^\eps \vert =O(1) \hbox{ for } i \in T_{i_0}\cup \{ i_0 \}  ,\\
& \vert \nabla \tilde{B}_{i}^\eps \vert =o(1) \hbox{ for } i\not\in T_{i_0}\cup \{ i_0 \}, \\
&\vert \nabla \tilde{R}^\eps \vert =O \left( 1 \right),
\end{split}
\ee
and
\be
\begin{split}
\Delta \tilde{B}_i^\eps &=o(1), \\
\Delta\tilde{R}^\eps &= o(1),
\end{split}
\ee
 and the relation of quasi-conformality
\be
\begin{split}
&\langle (\tilde{B}_i^\eps )_x , (\tilde{B}_i^\eps)_y \rangle = o(1),\\
& \langle( \tilde{B}_i^\eps  )_x ,  (\tilde{B}_i^\eps  )_x \rangle- \langle (\tilde{B}_i^\eps  )_y , ( \tilde{B}_i^\eps )_y \rangle = o(1 ).
\end{split}
\ee
and 
\be
\begin{split}
&\left\langle  \sum_{i} (\tilde{B}_i^\eps +  \tilde{R}^\eps )_x , \sum_{i} (\tilde{B}_i^\eps +  \tilde{R}^\eps )_y \right\rangle = o(1),\\
&\left\langle\sum_{i}( \tilde{B}_i^\eps +  \tilde{R}^\eps )_x ,  \sum_{i} (\tilde{B}_i^\eps +  \tilde{R}^\eps )_x \right\rangle- \left\langle  \sum_{i} (\tilde{B}_i^\eps +  \tilde{R}^\eps )_y , \sum_{i}( \tilde{B}_i^\eps +  \tilde{R}^\eps )_y \right\rangle = o(1 ).
\end{split}
\ee
Then, thanks to standard elliptic theory, see \cite{GT}, we get that   $\tilde{R}^\eps$ and $ \tilde{B}_{i}^\eps$ converge in $C^2_{loc} (\R^2\setminus\{ \tilde{S}_{i_0} \} )$ to   $\tilde{R}$ and $ \tilde{B}_{i}$ which satisfy
\be
\Delta \tilde{R}= \Delta \tilde{B}_{i}=0,
\ee
 the relations of conformality
\be
\begin{split}
&\langle (\tilde{B}_i )_x , (\tilde{B}_i)_y \rangle = 0,\\
& \langle( \tilde{B}_i)_x ,  (\tilde{B}_i )_x \rangle- \langle (\tilde{B}_i  )_y , ( \tilde{B}_i )_y \rangle = 0,
\end{split}
\ee
and 
\be
\begin{split}
&\left\langle  \sum_{i} (\tilde{B}_i +  \tilde{R})_x , \sum_{i} (\tilde{B}_i +  \tilde{R} )_y \right\rangle =0,\\
& \left\langle\sum_{i}( \tilde{B}_i +  \tilde{R} )_x ,  \sum_{i} (\tilde{B}_i +  \tilde{R} )_x \right\rangle- \left\langle  \sum_{i} (\tilde{B}_i +  \tilde{R} )_y , \sum_{i}( \tilde{B}_i +  \tilde{R} )_y \right\rangle = 0 .
\end{split}
\ee
On the one hand, if $i \in T_{i_0}\cup \{ i_0 \} $ we easily check that $\tilde{B}_i$ is a conformal parametrization of a plane. Then either this parametrization is singular if $\tilde{a}_{i} =\ds \lim_{\eps\rightarrow 0}\frac{ a_{i}^\eps -a_{i_0}^\eps}{d^\eps}$ is finite, or it is an affine map from $\R^2$ to $\R^3$.\\ 
On the other hand, thanks to the fact that  $r^\eps =O(t^\eps)$, then  $\nabla \tilde{R}$ is uniformly bounded on $\R^2$, hence  by the Liouville theorem $\nabla  \tilde{R}$ is constant, then $\tilde{R}$ is the standard parametrization of a plane. Let $j_0$ such that $\vert a_{i_0}^\eps -a_ {j_0}^\eps \vert =d^\eps$. Then we have the sum of at least two planes ($\tilde{B}_{i_0}$ and $\tilde{B}_{j_0}$) whose parametrization is singular in different points, which satisfies the equation of minimal surfaces. But, since this planes come from the limit of  embedded surfaces they must be parallel. Hence up to change the coordinates we can assume that the third coordinate of  $\nabla \left(\sum_{i} \tilde{B}_i +  \tilde{R} \right)$ vanishes. Then we have a conformal maps from $\R^2$ into itself with at least two singularities, we easily see that the minimal surface parametrized by $\sum_i \tilde{B}_i +  \tilde{R}$ necessary get a branched point. Here, the idea come from the Enneper-Weierstrass representation of minimal surfaces, see \cite{Osser}. Indeed, let $u$ be a solution of the minimal surface equation, we set 

$$\Phi = u_x + i u_y .$$
Then,  $\Phi$ is holomorphic and $\Phi^2 =0$, but applying this to $\tilde{B}_{i}$, we easily see that $\Phi$ is a rational fraction, with a pole if the parametrization is singular. Since we get at least two different poles, then this proves that $( \tilde{B}_i +  \tilde{R} )_x +i ( \tilde{B}_i +  \tilde{R} )_y$ vanishes some where. Finally the limit surface can be seen as a rational fraction of $\C$ whose derivative vanishes some where.\\ 

But, applying lemma \ref{simple}, we get a contradiction on the fact that $u^\eps$ is embedded, which achieves the proof of the claim.\hfill$\square$\\

%%%%%%%%%%%%%%%%%%%%%%%%%%%%%%%%%%%%%%%%%%%%%%%%%%%%%%%%%%%%%%%%%%%%%%  

%%%%%%%%%%%%%%%%%%%%%%%%%%%%%%%%%%%%%%%%%%%%%%%%%
%%%%%%%%%%%%%%%%%%%%%%%%%%%%%%%%%%%%%%
\appendix

\section{Why the bubbles are simple?}
\label{a1}
In this appendix, we give an explanation of the fact that an embedding can't converge to a branched surface.
\begin{lemma}
\label{simple}
Let $u^\eps :B(0,1) \rightarrow \R^3$ a sequence of smooth embedding such that there exists $u^0\in C^1(B(0,1),\R^3 )$ and
\be
u^\eps \rightarrow u^0 \hbox{ in } C^{2}_{loc} (B(0,1)\setminus \{ 0\} .
\ee

Then $u^0 $ can't be a multiple parametrization, that is to say there is no embedded $U_0\in C^1(B(0,1),\R^3 )$,  $\Phi \in \mathcal{O}(B(0,1),\C )$ an holomorphic function and an integer $k\geq 2$ such that
\be
\begin{split}
&u^0 = U^0 \circ \Phi\\
& \hbox{ and } \\
 &\Phi(z) = z^k+ o(\vert z\vert^k) \hbox{ as } z\rightarrow 0.
 \end{split}
 \ee
\end{lemma}

{\it Proof of the lemma \ref{simple} :} \\

First of all, up to a diffeomorphism of a neighborhood of $0$, we can assume that 
\be
u^\eps \rightarrow U^0(z^l) \hbox{ in } C^{2}_{loc}(B(0,\delta)\setminus\{0\}) .
\ee
where $l\geq 2$ and $\delta>0$. Let $A_\delta = B \left(0,\frac{\delta}{2}\right)\setminus B\left(0,\frac{\delta}{3}\right)$ and $C_r$ be the cylinder of center $U^0(0)$, radius $r$ and orthogonal to $T_{U^0(0)} U^0$, the tangent plane to the image of $U^0$ at $U^0(0)$. Let  $\delta>0$ and $r>0$ be small enough such that $C_{r}\cap U_0(A_\delta)$ is a simple curve. Then, for $\eps$ small enough, we easily see that the intersection of $u^\eps(A_\delta)$ and $C_r$ turn $l$ times around the cylinder, hence  $u^\eps(A_\delta)$ necessary intersect, which is a contradiction and proves the lemma.\hfill$\square$
%%%%%%%%%%%%%%%%%%%%%%%%%%%%%%%%%%%%%%%%%%%%%%%%%%%%%%
%%%%%%%%%%%%%%%%%%%%%%%%%%%%%%%%%%%%%%%%%%%
\section{Expansion of the metric and the Christoffel symbol}
\label{a2}
Using the classical expansion of the metric in a normal coordinates centered at $p\in N$, see \cite{Sakai96},  we get 
\be
g_{ij}(y)=\delta_{ij}+\frac{R_{ikmj}(p)}{3}y^{k}y^{m} +\frac{R_{ikmj,n}(p)}{6} y^{k}y^{m}y^{n}+o(r^3).
\ee
where $r^2 =\sum y_i^2$. Let $g^\eps(y)= g(\eps y)$ then we get
\be
(g_\eps)_{ij}(y)=\delta_{ij}+\frac{\eps^2}{3}R_{ikmj}(p) y^{k}y^{m} +\frac{\eps^3}{6} R_{ikmj,n}(p)  y^{k}y^{m}y^{n}+o(\eps^3).
\ee
Then we easily gets
\be
g_{\epsilon}^{ij} (y) =\delta_{ij}-\frac{\eps^2}{3}R_{ikmj}(p)y^{k}y^{m} -\frac{\eps^3}{6} R_{ikmj,n}(p)  y^{k}y^{m}y^{n}+o(\eps^3)
\ee
and
\be
\sqrt{\vert g_\eps \vert} (y) = 1- \frac{\eps^2}{6} \mathrm{Ric}_{mn}(p)y^m y^n - \frac{\eps^3}{12} \mathrm{Ric}_{mn,k}(p)y^m y^n y^k + o(\eps^3).
\ee
Now we are going to compute the expansion of the Christoffel symbol, using its expression with respect to  the metric, that is to say
\be
\Gamma_{ij}^{k}=\frac{1}{2}g^{kl}\left(g_{jl,i}+g_{il,j}-g_{ij,l}\right) .
\ee
Using the above formulas and the  second Bianchi identities, we get
\be
\begin{split}
(\Gamma_\eps)_{ij}^{k}(y) &= \frac{\eps^2}{3}\left( R_{jmik}(p) + R_{imjk}(p) \right) y^{m} \\
&+\frac{\eps^{3}}{6} \left(R_{jmik,n}(p) + R_{imjk,n}(p) +R_{imjn,k}(p) \right) y^{m}y^{n}+o(\eps^{3}),
\end{split}
\ee
where $(\Gamma_\eps)_{ij}^{k}$ is the Christoffel symbol of $g_\eps$. What we re-write it in more digest form   
\beq
\label{Gamma}
(\Gamma_\eps)_{ik}^{j} (y) = A_{ijkm}(p) y^{m}\eps^2 + B_{ijkmn}(p) y^{m}y^{n} \eps^3+o(\eps^{3})
\eeq
where $A_{ijkm}(p)=\frac{1}{3}\left( R_{kmij}(p) + R_{imkj}(p)\right)$ and\\ 
$B_{ijkmn}(p) = \frac{1}{12}\left( 2R_{kmij,n}(p)+ 2R_{imkj,n}(p) + R_{kmnj,i}(p) + R_{imnj,k}(p) - R_{imnk,j}(p) \right)$.

%%%%%%%%%%%%%%%%%%%%%%%%%%%%%%%%%%%%%%%%%%%%
\section{Linearized equation}
\label{a3}
Before to state our main result about the linearized equation, we give a lemma about the solution of $\Delta \alpha =\frac{8 }{(1+\vert x\vert^2)^2} \alpha $ which satisfy a decreasing assumption.
\begin{lemma}
\label{ChenLin}
Let $\alpha$ be a smooth solution of 
\beq
\label{lap2}
\begin{cases}
\Delta \alpha = \frac{8}{(1+\vert x\vert^2)^2}\alpha \hbox{ on } \R^2 ,\\
\alpha(0) =\nabla \alpha (0) =0 .
\end{cases}
\eeq
If $\vert \alpha (x) \vert \leq c(1+ \vert x \vert)^\tau$ for some $\tau\in [0,2[$ in $\R^2$, then $\alpha \equiv 0$.
\end{lemma}
Such a lemma has already be proved by Chen and Lin  with $\tau \in [0,1[$, see lemma 2.3 of \cite{ChenLin}. This result is not surprising, since the hypothesis of Chen and Lin  is almost equivalent to $\alpha  \in H^{1}(S^2)$ and in that case the main equation is just the  stereographic projection of $\Delta_{S^2} \alpha =2\alpha$ which has only zero as a solution which satisfies our initial data.\\

{\it  Proof of lemma \ref{ChenLin}:}\\

Here we repeat the proof of Chen and Lin with our additional estimate in order to get a larger set of admissible functions. In fact we prove first that the Fourier coefficient decrease more than expected.\\

Let $k\geq 2$ and  
$$ \alpha_k +i \beta_k = \int_{0}^{2\pi} \alpha e^{i k\theta} d\theta . $$
Then
\beq\label{alpha1}
\Delta \alpha_k = \left( \frac{8}{(1+r^2)^2} - \frac{k^2 }{r^2}\right)\alpha_k ,
\eeq
where $ \Delta \alpha_k  =\frac{-1}{r} \partial_r (r \partial_r \alpha_k)$. Then we set $\gamma_k =\frac{\alpha_k}{r^k}$ on $[1,+\infty[$, and we easily get that
\beq\label{alpha2}
\Delta \alpha_k = -k^2 r^{k-2} \gamma_k -(2k+1)r^{k-1}\gamma_k' -r^k \gamma_k'' .
\eeq
On the other hand, thanks to our hypothesis and (\ref{alpha1}), there exists $c$ a positive constant such that
\beq\label{alpha3}
\Delta \alpha_k \leq cr^{\tau-4} - k^2 r^{k-2}\gamma_k\ \hbox{ on } [1,+\infty[ .
\eeq
Hence, thanks to (\ref{alpha2}) and (\ref{alpha3}), we get 
\be
\begin{split}
  -(2k+1)r^{k-1}\gamma_k' -r^k \gamma_k'' &\leq c r^{\tau-4} \\
r^{2k+1}\gamma_k''  + (2k+1) r^{2k} \gamma_k' &\geq -c r^{\tau+k-3} \\
(r^{2k+1}\gamma_k')' &\geq -c r^{\tau+k-3} .
\end{split}
\ee
Then we integrate on $[1,r]\subset [1,+\infty[$, which gives
\be
\begin{split}
r^{2k+1}\gamma_k'(r) -\gamma_k'(1) &\geq c \frac{1-r^{\tau+k-2}}{\tau+k-2}\\
\gamma_k'(r) &\geq \left(\frac{1}{r}\right)^{2k+1}\gamma_k'(1) + \frac{c}{\tau+k-2} \left(\frac{1}{r^{2k+1}} - r^{\tau-k-3}\right) .
\end{split}
\ee
Then we integrate on $[R,+\infty[$, which gives
$$-\gamma_k(R) \geq \frac{R^{-2k}}{2k} \gamma_k'(1) + \frac{c}{\tau+k-2}\left( \frac{1}{(2k) R^{2k}} + \frac{R^{\tau-k-2}}{\tau-k-2}\right)$$

Here we used the fact that $\gamma_k(r)=O\left( r^{\tau-2}\right)$ and $\tau<2$. Thanks to last inequality, there exists $C$ a positive constant, such that  
$$\gamma_k(R) \leq C R^{\tau-k-2}$$
Then we get
$$\alpha_k(r) \leq C(r^{\tau-2}+1) \hbox{ on } [0,+\infty[.$$
Since the equation is linear we can applied the same argument to $-\alpha$ an finally we get the improved estimate, for every $k\geq 2$, there  exists $C_k$ a positive constant, such that  
$$\vert \alpha_k\vert  \leq C_k (1+r)^{\tau-2} \hbox{ on } [0,+\infty[ .$$
Of course the same result is true considering  $\beta_k$. Now we can follow the proof of Chen and Lin.\\

Let 
$$\psi_i (x) = \frac{x_i}{(1+\vert x\vert^2)} \hbox{ for } i=1,2 \hbox{ and } \psi_0(x)= \frac{1-\vert x\vert^2}{1+\vert x \vert^2}$$
 We are going to prove that any solution $\alpha$ of (\ref{lap2}) witch satisfies  $\vert \alpha (x) \vert \leq c(1+ \vert x \vert)^\tau$ for some $\tau\in ]0,2[$, is a linear combination of this three elementary solutions of \ref{lap2}, that is to say
 $$\alpha = \sum_{i=0}^{2} a_i \psi_i ,$$
 for some constant $a_i\in \R$. And then, the initial condition will give the result.\\
 
In order to show our result it suffices to show that $\alpha_k\equiv 0$ and $\beta_k\equiv 0$ for $k\geq 2$. We are going to prove this result for the $\alpha_k$, the argument are exactly the same for the $\beta_k$. Let $\phi_1 =  \int_{0}^{2\pi} \psi_1 cos(\theta) d\theta $. Then $\phi_1=O\left(\frac{1}{\vert x \vert}\right)$. Now, we suppose that $\alpha_k\not\equiv 0$ for some $k\geq 2$. Since $\phi_1(r) >0$ on $]0,+\infty[$, then by comparison with $\phi_1$, $\alpha_k$ never vanishes on $]0,+\infty[$.

Then, thanks to (\ref{alpha1}), we get  
\be
\begin{split}
\phi_1(r)\alpha_k'(r) r-\alpha_k(r)\phi_1'(r) r &= \int_{0}^{r}(\alpha_k\Delta\phi_1 - \phi_1 \Delta \alpha_k)s ds \\
&= (k^2-1) \int_{0}^{r} \frac{\alpha_k \phi_1}{s} ds .
\end{split}
\ee
Since $\vert \alpha_k \vert = O (1+r)^{\tau -2}$, then for a given positive content $C$, there exists a sequence of $r_i \rightarrow + \infty$ such that $\alpha_k'(r_i)r_i \leq C r_i^{\tau-2}$. Thus
$$
0= \lim_{i\rightarrow +\infty} \phi_1(r_i)\alpha_k'(r_i) r_i -\alpha_k(r_i)\phi_1'(r_i) r_i =(k^2-1) \int_{0}^{+\infty} \frac{\alpha_k \phi_1}{s} ds .
$$
Note that $\displaystyle  \frac{\alpha_k \phi_1}{s} =O((1+s)^{\tau-3})$ is integrable. Thus $\alpha_k\equiv 0$, which is a contradiction and prove the lemma.\hfill$\square$\\

%%%%%%%%%%%%%%%%%%%%%%%%%%%%%%%
Now we show how the study of the linearized problem can be reduce to the study of the previous equation.
\begin{prop}
\label{plin}
Let $\omega$ a simple solution of (\ref{eqlim}). Let also $r \in C^2(\R^2)$ be a solution of 
\beq
\label{eqlinea}
\Delta r +2 \left( r_ x \wedge \omega_y +\omega_ x \wedge r_ y \right)= 0 
\eeq
with
\beq
\label{linconf}
\begin{split}
& \langle r_x, \omega_x\rangle - \langle r_y , \omega_y\rangle=  0  ,\\
& \langle r_x, \omega_y\rangle +   \langle r_y, \omega_x\rangle= 0. 
\end{split}
\eeq
Setting $a,b,c,d,e$ and $f$ smooth functions as
\beq
\label{ws}
\nabla r =\left(\begin{array}{c}r_x \\r_y\end{array}\right) =\left(\begin{array}{c} a \omega_x +b \omega_y + c (\omega_x\wedge \omega_y )  \\ d \omega_x + e \omega_y +  f (\omega_x\wedge \omega_y) \end{array}\right),
\eeq
then $a,b,c,d,e$ and $f$ satisfy
\be
\begin{split}
&e= a\\
&d=-b\\
& \Delta a = \vert \nabla \omega \vert^2 a  \\
& \Delta b = \vert \nabla \omega \vert^2 b \\
& c= \frac{2}{\vert\nabla \omega\vert^2}\left( -a_x+b_y \right) \\
& \hbox{ and }\\
&f =  \frac{2}{\vert\nabla \omega\vert^2}\left(  -b_x- a_y  \right)
\end{split}
 \ee
\end{prop}

{\it Proof of proposition \ref{plin}:}\\

First of all, we easily get, thanks to (\ref{linconf}) and  (\ref{ws}), that
\be
\begin{split}
&a= e\\
&b=-d
\end{split}
 \ee

Differentiating (\ref{ws}), we get
\be
\begin{split}
\Delta r  &=\left( -a_x+b_y +c \frac{\vert \nabla \omega \vert^2}{2} \right) \omega_x +\left( -b_x- a_y  +f  \frac{\vert \nabla \omega \vert^2}{2} \right) \omega_y\\
&+ \left( -c_x- f_y  -2a - c \frac{(\vert \nabla \omega \vert^2)_x}{\vert \nabla \omega \vert^2} - f \frac{(\vert \nabla \omega \vert^2)_y}{\vert \nabla \omega \vert^2}  \right) (\omega_x\wedge \omega _y)\\
\end{split}
\ee
Now using (\ref{eqlinea}) and identifying  each coefficient of the equation in our special orthogonal frame, we get
\beq
\label{l1}
-a_x+ b_y  -c \frac{\vert \nabla  \omega \vert^2}{2}  = 0 ,
 \eeq
\beq
\label{l2} 
-b_x- a_y -f \frac{\vert \nabla \omega \vert^2}{2} = 0 .
\eeq
\beq
\label{l3}
-c_x- f_y +2a- c \frac{(\vert \nabla \omega \vert^2)_x}{\vert \nabla \omega \vert^2} - f \frac{(\vert \nabla \omega \vert^2)_y}{\vert \nabla \omega \vert^2} = 0
\eeq
Moreover, thanks to the fact that $r_{xy}= r_{yx}$ we get 
\beq
\label{l4}
a_y +b_x +f \frac{\vert \nabla \omega \vert^2}{2}  = 0 ,
 \eeq
\beq
\label{l5} 
b_y- a_x  -c \frac{\vert \nabla \omega \vert^2}{2} = 0 ,
\eeq
\beq
\label{l6}
2b +c_y- f_x  + c \frac{(\vert \nabla \omega \vert^2)_y}{\vert \nabla \omega \vert^2} - f \frac{(\vert \nabla \omega \vert^2)_x}{\vert \nabla \omega \vert^2} = 0 .
\eeq
Then, summing $(\ref{l1})_x$, $-(\ref{l4})_y$ and $- \frac{\vert \nabla \omega \vert^2}{2} (\ref{l3})$ we get
\be
\begin{split}
\Delta a - \vert \nabla \omega \vert^2 a  &=0 .
 \end{split}
 \ee
Then, summing $(\ref{l2})_x$, $-(\ref{l5})_y$ and $-\frac{\vert \nabla \omega \vert^2}{2} (\ref{l6})$ we get
\be
\begin{split}
\Delta b - \vert \nabla \omega \vert^2 b  &=0 .
 \end{split}
 \ee

Finally, thanks (\ref{l1}), (\ref{l2}) , we get 
\be
c= \frac{2}{\vert\nabla \omega\vert^2}\left( -a_x+b_y \right)
 \ee
 and
\be
f =  \frac{2}{\vert\nabla \omega\vert^2}\left(  -b_x- a_y  \right) .
\ee
\hfill$\square$\\

Here is our main result on the linearized problem. This classification is an improvement of existing results, see Lemma 9.1 of \cite{ChanilloMalchiodi} and Corollary 1.8 of \cite{CaldiroliMusina04}.
\begin{prop}
\label{lin}
let $\omega$ be a simple solution of (\ref{eqlim}) and  $r\in C^2(\R^2)$ be a solution of 
 \be
 \begin{split}
&\Delta r + 2 \left( r_ x \wedge \omega_y +\omega_ x \wedge r_ y \right)=0,\\
& \langle r_x, \omega_y\rangle +  \langle \omega_x, r_y\rangle  =0,\\
& \langle r_x, \omega_x\rangle =  \langle r_y , \omega_y\rangle  =0,\\
& \nabla r(0) =\nabla^2 r (\nabla \omega)(0) =0.
\end{split}
\ee
If $\vert \nabla r\vert $ is bounded  then $r$ is a constant function.
\end{prop}

{\it Proof  of proposition \ref{lin} :}\\

First of all, up to compose with an homography we can assume that 
$$ \omega(x,y)= \frac{1}{1+r^2}\left(\begin{array}{c}2x \\2y \\r^2-1\end{array}\right). $$
Indeed our equations are invariant with respect to a conformal transformation. Just the initial condition change in 
\be
 \nabla r(a) =\nabla^2 r(\nabla \omega)(a) =0.
\ee
where $a$ is the preimage of $0$  by the homography, of course it could be $\infty$. Now, up to compose by an  inversion or a translation, our initial condition comes back to zero.\\

We improve the decreasing assumption using Green function. Indeed Thanks to lemma \ref{green} and \ref{estimint1}, we easily get  that 
\beq
\label{estimgg}
\vert \nabla r (z)\vert = O\left(\frac{Ln(\vert z \vert)}{\vert z \vert}\right) \hbox{ when } z\rightarrow +\infty
\eeq
Now, let $a,b,c,d,e$ and $f$ as in the previous proposition. Then they satisfies
\be
\begin{split}
 &\Delta a = \vert \nabla \omega \vert^2 a  \\
& \Delta b = \vert \nabla \omega \vert^2 b \\
& c= \frac{2}{\vert\nabla \omega\vert^2}\left( -a_x+b_y \right) \\
& \hbox{ and }\\
&f =  \frac{2}{\vert\nabla \omega\vert^2}\left(  -b_x- a_y  \right)\\
&e-a=d+b=0 .
\end{split}
 \ee
Moreover thanks to our initial condition and (\ref{estimgg}), we have 
\be 
\begin{split} 
& a(0)=b(0)=\nabla a(0)=\nabla b(0)=0 , \\
& \vert a \vert , \vert b \vert =O( (1 + \vert z \vert )^\frac{3}{2} ). 
\end{split}
\ee
Then $a$ and $b$ satisfies the hypothesis of lemma \ref{ChenLin}, then $a\equiv b\equiv 0$ and we easily prove that $r$ is a constant function . \hfill$\square$\\
%%%%%%%%%%%%%%%%%%%%%%%%%%%%%%%%%%%%%%%%%%%%

%%%%%%%%%%%%%%%%%%%%%%%%%%%%%%%%%%%

 \section{Green functions and integral estimates}
 \label{a4}
 
 Let $G$ and $G_R$ be the Green function of the Laplacian respectively  on the plane and on the ball $B(0, R)$, that is to say
 \be
 \begin{split}
 & G(z_1,z_2) = \frac{1}{2\pi} Ln \vert z_1 - z_2\vert ,\\
 &G_R(z_1,z_2) = \frac{1}{2\pi} \left( Ln \vert z_1 - z_2\vert  -Ln \left\vert \frac{R}{\vert z_1\vert}z_1  -\frac{\vert z_1\vert}{R} z_2\right\vert \right) .
 \end{split}
 \ee
  \begin{lemma}
 \label{green} 
 Let $u$ and $f$ two functions in $C^2(\R^2,\R)$ which satisfies
 \be
 \begin{cases}
 \Delta u = f,\\
 \Vert \nabla u \Vert_{\infty} <+\infty ,\\
 f= O\left(\frac{1}{\vert z \vert^2 }\right) .
 \end{cases} 
 \ee
 Then we have
 $$ \nabla u(z_0) = \int_{\R^2} \nabla G(z_0,z) f(z) dz .$$ 
 \end{lemma}

 {\it Proof of lemma \ref{green}:}\\
 
Let $z_0 \in \R^2$ and $R>0$ such that  $z_0 \in B(0,R)$, then thanks to the standard  Green formula, we get 
 
$$ \nabla u(z_0) = \int_{B(0, R)} G_R(z_0,z) \nabla f(z) dz + \int_{\partial B(0, R)} \frac{\partial  G_R}{\partial n}(z_0,z) \nabla u(z) d\sigma .$$
 Then, integrating by part, we get  
  \beq
 \label{gg}
 \begin{split}
  \nabla u(z_0) &= \int_{B(0, R)} \nabla G_R(z_0,z) f(z) dz  +  \int_{\partial B(0, R)} G_R(z_0,z) f(z) dz \\
   &+ \int_{\partial B(0, R)} \frac{ \partial G_R}{\partial n}(z_0,z) \nabla u(z) d\sigma .
   \end{split}
  \eeq
 For $z_0$ fixed, we get
\be
\begin{array}{cc} \begin{cases}
\vert  G_R(z_0,z) \vert = O\left(Ln\vert z - z_0 \vert \right) \\
\vert \nabla G_R(z_0,z) \vert = O\left(\frac{1}{\vert z - z_0 \vert}\right)
\end{cases} & \hbox{ when } z\rightarrow +\infty,  \end{array}
\ee
  and
 $$ \left\vert \nabla \left(\frac{G_R}{\partial n}\right)(z_0,z) \right\vert= O\left(\frac{1}{R^2}\right) \hbox{ when } z\in \partial B(0,R) \hbox{ and } R\rightarrow +\infty  .  $$ 
 
This  allows us  to take the limit when  $R$ goes to infinity in (\ref{gg}) and this gives the result.\hfill$\square$\\
 
 \begin{lemma}
 \label{estimint1} There exits a positive constant $C$ such that, for all $z_0\in \R^2$,
$$\int_{\R^2} \vert \nabla G(z,z_0) \vert \frac{1}{1+\vert z \vert^2} dz \leq C \frac{Ln(2+\vert z_0\vert)}{1+\vert z_0\vert} .$$
\end{lemma}

{\it Proof of lemma \ref{estimint1} :}\\

Applying again standard estimates on Green functions, there exists  a positive constant $C$, such that
\be
\begin{split}
\int_{\R^2} \vert \nabla G(z,z_0)\vert  \frac{1}{1+\vert z \vert^2} dz & \leq \int_{\R^2} \frac{C}{\vert z-z_0\vert(1+\vert z \vert^2)}  dz \\
& =\int_{B\left(z_0, \frac{\vert z_0\vert}{2}\right)} \frac{C}{\vert z-z_0\vert(1+\vert z \vert^2)}  dz \\
&+ \int_{B\left(0, \frac{\vert z_0\vert}{2}\right)} \frac{C}{\vert z-z_0\vert(1+\vert z \vert^2)}  dz \\
&+ \int_{\left\{ \vert z \vert \geq \frac{\vert z_0\vert}{2},\vert z-z_0 \vert \geq \frac{\vert z_0\vert}{2}\right\}} \frac{C}{\vert z-z_0\vert(1+\vert z \vert^2)} dz \\
& \leq \frac{4C}{4+\vert z_0\vert^2} \int_{B\left(z_0, \frac{\vert z_0\vert}{2}\right)} \frac{1}{\vert z-z_0\vert}  dz \\
&+ \frac{C}{\vert z_0 \vert} \int_0^{\frac{\vert z_0\vert}{2}} \frac{2r}{1+r^2} dr\\
&+ 4C \int_{\frac{\vert z_0\vert}{2}}^\infty \frac{1}{r^2} dr\\
& \leq C\left( \frac{2\vert z_0\vert}{4+\vert z_0\vert^2} +  \frac{1}{\vert z_0 \vert} Ln \left( 1+\frac{\vert z_0 \vert^2 }{4}\right) +   \frac{8}{\vert z_0 \vert}\right), 
\end{split}
\ee
which proves the lemma.\hfill$\square$\\
%%%%%%%%%%%%%%%%%%%%%%%%%%%%%%%%%%%%%%%%%%%%
\section{Wente inequality and applications} 
\label{a5}

First of all, we remind us the Wente inequality. Thanks to the work of Bethulel, Ghigladia  and Topping, see \cite{BethuelGhidaglia} and \cite{Topping97}, we have the following version of the Wente inequality.
\begin{thm}[Wente inequality]
\label{wente1}
Let $\Omega$ be a bounded open set of $\R^2$ and $v\in H^1(\Omega)$. Let $u\in W^{1,1}_0(\Omega)$ be the solution of
$$\Delta u = -2  v_x \wedge v_y  \hbox{ on } \Omega ,$$
then 
$$\Vert u\Vert_{\infty} +\Vert \nabla u \Vert_2 \leq \frac{1}{\pi} \Vert \nabla v\Vert_2^2 .$$ 
\end{thm}
Which is remarkable here is that the constant is independent of  $\Omega$. Then using such a result, Topping as proved a Wente's inequality for surfaces, see theorem 4 of \cite{Topping97}.
\begin{thm}
\label{topp} Let $\Sigma$ a compact Riemannian surface and $v\in H^1(\Sigma, \R^2)$. Then if $u\in W^{1,1}(\Sigma)$ be the solution of

$$\Delta u = det(\nabla v) \hbox{ on } \Sigma ,$$
then 
$$osc(u)+\Vert \nabla u \Vert_2 \leq \frac{1}{\pi} \Vert \nabla v\Vert_2^2 ,$$
where $\ds osc(u)= \sup_{x,y\in \Sigma}\vert u(x)-u(y)\vert .$ 
\end{thm}
Then,  assuming that $u\in H^1$, we extend such an equality to $\Omega=\R^2$.
\begin{cor}
\label{wente2}
Let $v\in H^1(\R^2)$ and  $u\in H^{1}(\R^2)$ be a solution of
$$\Delta u = -2  v_x \wedge v_y \hbox{ on } \R^2$$
then
$$\Vert \nabla u \Vert_2 \leq \frac{2}{\pi}  \Vert \nabla v\Vert_2^2 .$$ 
\end{cor}

{\it Proof of corollary \ref{wente2} :}\\

Let $\pi$ the standard stereographic projection from $S^2$ to $\R^2$. Thanks to the conformal invariance of the equation, $u \circ \pi^{-1}$ and $v\circ \pi^{-1}$ satisfies the hypothesis of  theorem \ref{topp} when $\Sigma=S^2$, hence we get that 
$$osc(u) \leq  \frac{1}{\pi}  \Vert \nabla v\Vert_2^2 $$
Then testing the equation against $u$ and integrating by parts we get the desired inequality.\hfill$\square$\\

Now we are in position to prove the our lemma, which allowed to control the supremum of the gradient for such solution. This is a new manifestation of the presence of a strong compensation phenomena in this equation.
\begin{lemma}
\label{mll}
Let $v\in H^{1}(\R^2)$ and  $u\in H^{1}(\R^2)$ be a solution of 
\beq
\label{ml}
\Delta u = v_x \wedge v_y. 
\eeq
Then there exists a positive constant $C$, independent of $v$, such that
$$\Vert \nabla u \Vert_\infty  \leq C \Vert \nabla v\Vert_\infty \Vert \nabla v \Vert_2 .$$ 
\end{lemma}

The proof of this lemma relies on the corollary \ref{wente2} and the following interpolation inequality.

\begin{lemma}[lemma A.2 \cite{BBH}]
Let $\Omega$ a smooth domain. Assume u satisfies 
$$ 
\begin{cases}
& \Delta u=f \hbox{ on } \Omega \\
& u  =  0  \hbox{ on } \partial\Omega 
\end{cases} 
$$
Then 
$$\Vert  \nabla u \Vert_{\infty}^2 \leq C \Vert f\Vert_{\infty}\, \Vert u\Vert_\infty $$ 
where $C$  is a  constant  depending  only  on  $\Omega$. 
\end{lemma}

To conclude this appendix, we give an other useful version of the Wente's inequality, see \cite{BC2} for example.

\begin{lemma}
\label{wente3}
Let $\Omega=B(0,1)$, $u\in H^1(\Omega)\cap L^\infty(\Omega) $ and $v\in H^1_0(\Omega)$, then there exists $C$, independent of $u$ and $v$, such that 
$$\left\vert \int_{\Omega} \langle u, v_x \wedge v_y \rangle \right\vert\leq C \Vert \nabla v\Vert_2  \Vert \nabla u\Vert_2^2 .$$
\end{lemma}

 \bibliographystyle{plain}
\bibliography{cmcbiblio}
\end{document}